\documentclass[11pt,reqno]{amsproc}
\usepackage{fullpage}
\numberwithin{equation}{section}
\usepackage{cite}
\usepackage{color}
\usepackage{graphicx}
\usepackage{subfig}
\usepackage{wrap fig,lip sum,booktabs}
\usepackage{caption}
\usepackage{comment}
\usepackage{algpseudocode}
\usepackage{algorithm}
\usepackage[scale=0.8]{geometry}
\usepackage{epstopdf}
\usepackage{enumerate}
\usepackage[debug=false, colorlinks=true, pdfstartview=FitV, linkcolor=blue, citecolor=blue, urlcolor=blue]{hyperref}
\newtheorem{remark}{Remark}

\newlength{\drop}
\definecolor{amethyst}{rgb}{0.6, 0.4, 0.8}
\definecolor{burgundy}{rgb}{0.5, 0.0, 0.13}

\title[Large-scale Non-negative Computational Framework]{Large-scale 
  Optimization-based Non-negative Computational Framework for 
  Diffusion Equations:~Parallel Implementation and Performance 
  Studies}

\author{\textbf{J.~Chang, S.~Karra and K.~B.~Nakshatrala} \\
  {\small Correspondence to: \textbf{\emph{e-mail:}} knakshatrala@uh.edu, 
  \textbf{\emph{phone:}}+1-713-743-4418}}

\date{\today}

\begin{document}

%================================================================================;
%  Below is CAML logo. All the students working in CAML will use the CAML logo.  ;
%  Do not change the format of the logo paper.                                   ;
%================================================================================;

\begin{titlepage}
    \drop=0.1\textheight
    \centering
    \vspace*{\baselineskip}
    \rule{\textwidth}{1.6pt}\vspace*{-\baselineskip}\vspace*{2pt}
    \rule{\textwidth}{0.4pt}\\[\baselineskip]
    {\LARGE \textbf{\color{burgundy} 
    Large-scale optimization-based non-negative \\[\baselineskip]
    computational framework for diffusion equations: \\[\baselineskip]
    Parallel implementation and performance studies}}\\[0.3\baselineskip]
    \rule{\textwidth}{0.4pt}\vspace*{-\baselineskip}\vspace{3.2pt}
    \rule{\textwidth}{1.6pt}\\[\baselineskip]
    \scshape
    An e-print of the paper is available on arXiv:~1506.08435. \par
    
    \vspace*{0.5\baselineskip}
    Authored by \\[0.5\baselineskip]
    
    {\Large J.~Chang\par}
    {\itshape Graduate Student, University of Houston.}\\[0.5\baselineskip]
    
    {\Large S.~Karra\par}
    {\itshape Staff Scientist, Los Alamos National Laboratory.}\\[0.5\baselineskip]
    
    {\Large K.~B.~Nakshatrala\par}
    {\itshape Department of Civil \& Environmental Engineering \\
      University of Houston, Texas 77204--4003.\\
    \textbf{phone:} +1-713-743-4418, \textbf{e-mail:} knakshatrala@uh.edu \\
    \textbf{website:} http://www.cive.uh.edu/faculty/nakshatrala\par}
    \vspace*{0.5\baselineskip}
\begin{figure}[h]
    \includegraphics[scale=0.29,clip]{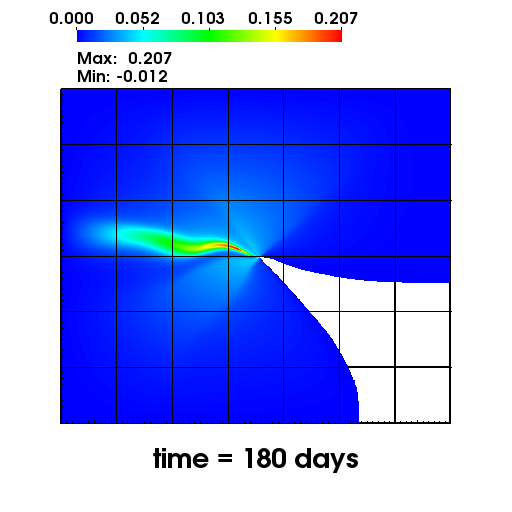}
    \caption*{\emph{This figure shows the fate of chromium after 180 days 
    using the single-field Galerkin formulation. The white regions indicate 
    the violation of the non-negative constraint.}}
\end{figure}
    \vfill
    {\scshape 2015} \\
    {\small Computational \& Applied Mechanics Laboratory} \par
  \end{titlepage}

%=========================;
%  Abstract and Keywords  ;
%=========================;
\begin{abstract}
It is well-known that the standard Galerkin formulation, 
which is often the formulation of choice under the finite 
element method for solving self-adjoint diffusion equations, 
does not meet maximum principles and the non-negative 
constraint for anisotropic diffusion equations. Recently, 
optimization-based methodologies that satisfy maximum 
principles and the non-negative constraint for steady-state 
and transient diffusion-type equations have been proposed. 
To date, these methodologies have been tested only on small-scale 
academic problems. The purpose of this paper is to systematically 
study the performance of the non-negative methodology in the 
context of high performance computing (HPC). PETSc and TAO libraries 
are, respectively, used for the parallel environment and optimization 
solvers. For large-scale problems, it is important for computational 
scientists to understand the computational performance of current algorithms 
available in these scientific libraries. 
The numerical experiments are conducted on the 
state-of-the-art HPC systems, and a single-core performance 
model is used to better characterize the efficiency of the solvers. 
Our studies indicate that the proposed non-negative computational 
framework for diffusion-type equations exhibits excellent strong 
scaling for real-world large-scale problems.
\end{abstract}
\maketitle

%========================;
%  Include all sections  ;  
%========================;

%*********************************************;
%                                             ;
%  NAME                                       ;
%  S1_LargeScale_Intro.tex                    ;
%                                             ;
%  WRITTEN BY                                 ;
%    Justin Chang                             ;
%    Satish Karra	         	      ;
%    Kalyana Babu Nakshatrala                 ;
%                                             ;
%*********************************************;
\section{INTRODUCTION}
\label{Sec:Introduction}
The modeling of flow and transport in subsurface is vital 
for energy, climate and environmental applications. Examples 
include CO$_2$ migration in carbon-dioxide sequestration, 
enhanced geothermal systems, oil and gas production, 
radio-nuclide transport in a nuclear waste repository, 
groundwater contamination, and  thermo-hydrology in the 
Arctic permafrost due to the recent climate change 
\cite{karra2014three,lichtner2014modeling,
hammond2010field,kelkar2014simulator}. 
Several numerical codes (e.g., FEHM \cite{zyvoloski2007fehm}, 
TOUGH \cite{pruess2004tough}, PFLOTRAN \cite{lichtner2015pflotran}) 
have been developed to model flow and transport in subsurface 
at reservoir-scale. These codes typically solve unsteady Darcy 
equations for flow and advection-diffusion equation for 
transport. The predictive capability of a numerical simulator 
depends on the robustness of the underlying numerical methods. 
A necessary and essential requirement is to satisfy important 
mathematical principles and physical constraints. 
One such property in transport and reactive-transport problems 
is that the concentration of a chemical species cannot be negative. 
Mathematically, this translates to the satisfaction of the 
discrete maximum principle (DMP) for diffusion-type equations. 
Subsurface flow and transport applications typically encounter 
geological media that are highly heterogeneous and anisotropic 
in nature, and it is well-known that the classical finite 
element (or finite volume and finite difference, for that 
matter) formulations do not produce non-negative solutions 
on arbitrary meshes for such porous media 
\cite{Ciarlet_Raviart_CMAME_1973_v2_p17,
liska2008enforcing,
Nakshatrala_JCP_2009_v228_p6726,
Nagarajan_IJNMF_2011_v67_p820}.  

%======================================================;
%  Subsection: Literature review on maximum principles ;
%======================================================;
Several studies over the years have focused on the 
development of methodologies that enforce the DMP 
and ensure non-negative solutions \cite{Nakshatrala_JCP_2009_v228_p6726,
Nagarajan_IJNMF_2011_v67_p820,Payette_IJNME_2012_v91_p742,
Nakshatrala_JCP_2013_v253_p278}. However, these studies 
did not address how these methods can be used for 
realistic large-scale subsurface problems that have millions of grid nodes. 
Furthermore, complex coupling between different physical 
processes as well as the presence of multiple species amplify the 
degrees-of-freedom (i.e., the number of unknowns).
The aim of this paper is to develop a parallel 
computational framework that solves anisotropic diffusion 
equations on general meshes, ensures non-negative solutions, 
and can be employed to solve large-scale realistic problems.

Large-scale problems can be tackled by using recent 
advancements in high-performance computing (HPC) 
methods and toolkits that can be used on the 
state-of-the-art supercomputing architecture. 
One such toolkit is PETSc \cite{petsc-user-ref}, which 
provides data structures and subroutines for setting up 
structured and unstructured grids, parallel communication, 
linear and non-linear solvers, and parallel I/O. These 
high-level data structures and subroutines help in faster 
development of parallel application codes and minimize 
the need to program low-level message passing, so that 
the domain scientists can focus more on the application. 
To this end, we develop a non-negative parallel framework
by leveraging the existing capabilities within PETSc. Our 
framework ensures the DMP for anisotropic diffusion by using 
lower-order finite elements and the optimization-based 
approach in \cite{Nakshatrala_JCP_2009_v228_p6726,
liska2008enforcing,Nagarajan_IJNMF_2011_v67_p820}. The TAO 
toolkit \cite{tao-user-ref}, which is built on top of PETSc, 
is used for solving the resulting optimization problems. The 
robustness of the proposed framework will be demonstrated by 
solving realistic large-scale problems. 

The rest of this paper is organized as follows. In Section 
\ref{Sec:Governing_Equations}, we present the governing 
equations and the classical single-field Galerkin finite 
element formulation for steady-state and transient diffusion 
equations. The optimization-based method to ensure non-negative 
concentrations is also outlined in this section. 
In Section \ref{Sec:Parallel_Implementation}, the parallel 
implementation procedure using PETSc and TAO is presented. 
We also highlight the relevant data structures used in 
this study and present a pseudo algorithm describing our parallel framework. 
In Section \ref{Sec:Performance_Modeling}, a performance model 
loosely based on the roofline model is outlined. This model is used 
to estimated the efficiency with respect to computing hardware 
of currently available solvers within PETSc and TAO.
In Section \ref{Sec:Numerical_Results}, we first verify our 
implementation using a 3D benchmark problem from the literature
and present a detailed performance study using the proposed model. 
Then, we study a large-scale three-dimensional realistic 
problem involving the transport of chromium in the subsurface
and document the numerical results of the non-negative 
methodology with the classical single-field Galerkin formulation. 
Conclusions are drawn in Section \ref{Sec:Concluding_Remarks}.

%*********************************************;
%                                             ;
%  NAME                                       ;
%    S2_LargeScale_GE.tex          	      ;
%                                             ;
%  WRITTEN BY                                 ;
%    Kalyana Babu Nakshatrala                 ;
%                                             ;
%*********************************************;
\section{GOVERNING EQUATIONS AND ASSOCIATED NON-NEGATIVE NUMERICAL METHODOLOGIES}
\label{Sec:Governing_Equations}
Let $\Omega \subset \mathbb{R}^{nd}$ be a bounded open 
domain, where ``$nd$'' is the number of spatial dimensions. 
The boundary of the domain is denoted by $\partial \Omega 
= \overline{\Omega} - \Omega$, which is assumed to be 
piecewise smooth. A spatial point is denoted by $\mathbf{x} 
\in \overline{\Omega}$. The gradient and divergence operators 
with respect to $\mathbf{x}$ are, respectively, denoted as 
$\mathrm{grad}[\cdot]$ and $\mathrm{div}[\cdot]$. As usual, 
the boundary is divided into two parts: $\Gamma^{\mathrm{D}}$ 
and $\Gamma^{\mathrm{N}}$. $\Gamma^{\mathrm{D}}$ is that part 
of the boundary on which Dirichlet boundary conditions 
are  prescribed, and $\Gamma^{\mathrm{N}}$ is the part of 
the boundary on which Neumann boundary conditions are 
prescribed. For mathematical well-posedness, we assume 
$\Gamma^{\mathrm{D}} \cup \Gamma^{\mathrm{N}} = \partial 
\Omega$ and $\Gamma^{\mathrm{D}} \cap \Gamma^{\mathrm{N}} 
= \emptyset$. 
The unit outward normal to boundary is denoted as 
$\widehat{\mathbf{n}}(\mathbf{x})$. The diffusivity 
tensor is denoted by $\mathbf{D}(\mathbf{x})$, which 
is assumed to symmetric, bounded above and uniformly 
elliptic. That is,
%------------------------------------------;
%  Equation: Symmetry of diffusion tensor  ;
%------------------------------------------;
\begin{align}
  \mathbf{D}(\mathbf{x}) = \mathbf{D}^{\mathrm{T}}
  (\mathbf{x}) \quad \forall \mathbf{x} \in \Omega
\end{align}
and there exists two constants $\mathrm{0 < 
  \xi_1 \leq \xi_2 < +\infty}$ such that 
%----------------------------------------;
%  Equation: Positive definiteness of D  ;
%----------------------------------------;
\begin{align}
  \label{Eqn:Helmholtz_positive_definiteness_D}
  \xi_1 \mathbf{y}^{\mathrm{T}} \mathbf{y}\leq \mathbf{y}^{\mathrm{T}} 
  \mathbf{D}(\mathbf{x}) \mathbf{y} \leq \xi_2 \mathbf{y}^{\mathrm{T}} 
  \mathbf{y} \quad \forall \mathbf{x} \in \Omega \; \mathrm{and} \; 
  \forall \mathbf{y} \in \mathbb{R}^{nd}
\end{align} 

%=============================================================;
%  Subsection: Governing equations for steady-state response  ;
%=============================================================;
\subsection{Governing equations for steady-state response}
We shall denote the steady-state concentration field 
by $c(\mathbf{x})$. The governing equations can be 
written as follows: 
%------------------------------------------------;
%  Equation: Steady-state anisotropic diffusion  ;
%------------------------------------------------;
\begin{subequations}
  \label{Eqn:LargeScale_BVP}
  \begin{alignat}{2}
  \label{Eqn:LargeScale_GE}
  -&\mathrm{div}[\mathbf{D}(\mathbf{x}) 
    \mathrm{grad}[c]] = f(\mathbf{x}) & \qquad & 
  \mathrm{in} \; \Omega  \\
  \label{Eqn:LargeScale_Dirichlet}
  &c(\mathbf{x}) = c^{\mathrm{p}} (\mathbf{x}) & \qquad &
  \mathrm{on} \; \Gamma^{\mathrm{D}} \\
  \label{Eqn:LargeScale_Neumann}
  -&\mathbf{\widehat{n}}(\mathbf{x}) \cdot \mathbf{D}
  (\mathbf{x}) \mathrm{grad}[c] =
  q^{\mathrm{p}}(\mathbf{x}) & \quad & \mathrm{on} \;
  \Gamma^{\mathrm{N}}
  \end{alignat}
\end{subequations}
where $f(\mathbf{x})$ is the volumetric 
source/sink, $c^{\mathrm{p}}(\mathbf{x})$ is the 
prescribed concentration, and $q^{\mathrm{p}}
(\mathbf{x})$ is the prescribed flux. For 
uniqueness, we assume $\Gamma^{\mathrm{D}} 
\neq \emptyset$. 

%====================================================================;
%  Subsubsection: Maximum principle and the non-negative constraint  ;
%====================================================================;
\subsubsection{Maximum principle and the non-negative constraint}
The above boundary value problem is a self-adjoint 
second-order elliptic partial differential equation 
(PDE). It is well-known that such PDEs possess an 
important mathematical property -- the classical 
maximum principle \cite{Evans_PDE}. The mathematical 
statement of the classical maximum principle can be 
written as follows: 
If $c(\mathbf{x}) \in C^{2}(\Omega) \cap C^{0}(\overline{\Omega})$, 
$\partial \Omega = \Gamma^{\mathrm{D}}$, and $f(\mathbf{x}) \leq 0$  
in $\Omega$ then 
\begin{align}
  \max_{\mathbf{x} \in \overline{\Omega}} c(\mathbf{x}) 
  = \max_{\mathbf{x} \in \partial \Omega} c^{\mathrm{p}}(\mathbf{x})
\end{align}
Similarly, if $f(\mathbf{x}) \geq 0$ in $\Omega$ then
%-------------------------------;
%  Equation: Maximum principle  ;
%-------------------------------;
\begin{align}
  \label{Eqn:LargeScale_MP}
  \min_{\mathbf{x} \in \overline{\Omega}} c(\mathbf{x}) 
  = \min_{\mathbf{x} \in \partial \Omega} c^{\mathrm{p}}(\mathbf{x})
\end{align}
To make our presentation on maximum principles 
simple, we have assumed stronger regularity on 
the solution (i.e., $c(\mathbf{x}) \in C^{2} 
\cap C^{0}(\overline{\Omega})$), and assumed that 
Dirichlet boundary conditions are prescribed on 
the entire boundary. However, maximum principles 
requiring milder regularity conditions on the 
solution, even for the case when Neumann boundary 
conditions are prescribed on the boundary, can be 
found in  literature (see 
\cite{Mudunuru_Nakshatrala_arXiv_2015,
Mudunuru_Nakshatrala_ADR_arXiv_2015}). 

If $f(\mathbf{x}) \geq 0$ in $\Omega$ and 
$c^{\mathrm{p}}(\mathbf{x}) \geq 0$ on the 
entire $\partial \Omega$ then the maximum 
principle implies that $c(\mathbf{x}) 
\geq 0$ in the entire domain, which is 
the non-negativity of the concentration 
field.  

%=========================================================;
%  Subsubsection: Single-field Galerkin weak formulation  ;
%=========================================================;
\subsubsection{Single-field Galerkin weak formulation}
The following function spaces will be 
used in the rest of this paper: 
%-----------------------------;
%  Equation: Function spaces  ;
%-----------------------------;
\begin{align}
  \mathcal{U} &:= \left\{c(\mathbf{x}) \in H^{1}(\Omega) 
  \; \big| \; c(\mathbf{x}) = c^{\mathrm{p}}(\mathbf{x}) \; 
  \mathrm{on} \; \Gamma^{\mathrm{D}}\right\} \\
  \mathcal{W} &:= \left\{w(\mathbf{x}) \in H^{1}(\Omega) 
  \; \big| \; w(\mathbf{x}) = 0 \; \mathrm{on} \; 
  \Gamma^{\mathrm{D}}\right\} 
\end{align}
where $H^{1}(\Omega)$ is a standard Sobolev space 
\cite{adams2003sobolev}. The single-field Galerkin 
weak formulation corresponding to equations 
\eqref{Eqn:LargeScale_GE}--\eqref{Eqn:LargeScale_Neumann} 
reads: Find $c(\mathbf{x}) \in \mathcal{U}$ 
such that we have 
%--------------------------------------;
%  Equation: Single field formulation  ;
%--------------------------------------;
\begin{align}
  \label{Eqn:Single_field_formulation}
  \mathcal{B}(w;c) = L(w) \quad \forall 
  w(\mathbf{x}) \in \mathcal{W}
\end{align}
where the bilinear form and linear functional 
are, respectively, defined as 
%-------------------------------------------------;
%  Equation: Bilinear form and linear functional  ;
%-------------------------------------------------;
\begin{subequations}
  \label{Eqn:functionals_B_L}
  \begin{align}
    \mathcal{B}(w;c) &:= \int_{\Omega} \mathrm{grad}[w(\mathbf{x})] 
    \cdot \mathbf{D}(\mathbf{x}) \mathrm{grad}[c(\mathbf{x})] 
    \; \mathrm{d} \Omega 
    \\
    L(w) &:= \int_{\Omega} w(\mathbf{x}) f(\mathbf{x}) 
    \; \mathrm{d} \Omega 
    + \int_{\Gamma^{\mathrm{N}}} w(\mathbf{x}) 
    q^{\mathrm{p}}(\mathbf{x}) \; \mathrm{d} \Gamma 
  \end{align}
\end{subequations}
Since $\mathbf{D}(\mathbf{x})$ is symmetric, 
by Vainberg's theorem \cite{Hjelmstad}, the 
single-field Galerkin weak formulation given 
by equation \eqref{Eqn:Single_field_formulation} 
is equivalent to the following variational 
problem: 
%-----------------------------------;
%  Equation: Variational statement  ;
%-----------------------------------; 
\begin{align}
  \label{Eqn:LargeScale_var_stat}
  \mathop{\mathrm{minimize}}_{c(\mathbf{x}) \in \mathcal{U}} 
  \quad \frac{1}{2} \mathcal{B}(c;c) - L(c)
\end{align}   

%================================================================;
%  Subsubsection: A methodology to enforce the maximum principle  ;
%================================================================;
\subsubsection{A methodology to enforce the maximum principle 
  for steady-state problems}
Our methodology is based on the finite element method. 
We decompose the domain into ``$Nele$'' non-overlapping 
open element sub-domains such that 
%-----------------------------------------;
%  Equation: Decomposition into elements  ;
%-----------------------------------------;
\begin{align}
  \overline{\Omega} = \bigcup_{e = 1}^{Nele} \overline{\Omega}^{e}
\end{align}
(Recall that a superposed bar denotes the set closure.) The 
boundary of $\Omega^e$ is denoted by $\partial \Omega^{e} 
:= \overline{\Omega}^{e} - \Omega^{e}$. We shall define the 
following finite dimensional vector spaces of $\mathcal{U}$ 
and $\mathcal{W}$:
%------------------------------------------;
%  Equation: Finite dimensional subspaces  ;
%------------------------------------------;
\begin{subequations}
  \begin{align}
    \mathcal{U}^{h} &:= \left\{c^{h}(\mathbf{x}) \in \mathcal{U} 
    \; \big| \; c^{h}(\mathbf{x}) \in C^{0}(\overline{\Omega}), 
    c^{h}(\mathbf{x}) \big|_{\Omega^e} \in \mathbb{P}^{k}(\Omega^{e}), 
      e = 1, \cdots, Nele \right\} \\
    \mathcal{W}^{h} &:= \left\{w^{h}(\mathbf{x}) \in \mathcal{W} 
    \; \big| \; w^{h}(\mathbf{x}) \in C^{0}(\overline{\Omega}), 
    w^{h}(\mathbf{x}) \big|_{\Omega^e} \in \mathbb{P}^{k}(\Omega^{e}), 
    e = 1, \cdots, Nele \right\}       
  \end{align}
\end{subequations}
where $k$ is a non-negative integer, and 
$\mathbb{P}^{k}(\Omega^{e})$ denotes the linear 
vector space spanned by polynomials up to $k$-th 
order defined on the sub-domain $\Omega^{e}$. The 
finite element formulation for equation 
\eqref{Eqn:Single_field_formulation} can be 
written as: Find $c^{h}(\mathbf{x}) \in \mathcal{P}^{h}$ 
such that we have 
%-----------------------------------------------------;
%  Equation: Single field finite element formulation  ;
%-----------------------------------------------------;
\begin{align}
  \label{Eqn:FE_formulation}
  \mathcal{B}(q^{h};c^{h}) = L(q^{h}) \quad 
  \forall q^{h}(\mathbf{x}) \in \mathcal{Q}^{h}
\end{align}
It has been documented in the literature that the 
above finite element formulation violates the 
maximum principle and the non-negative constraint \cite{liska2008enforcing,
Nakshatrala_JCP_2009_v228_p6726,Nagarajan_IJNMF_2011_v67_p820}.
 
We now outline an optimization-based methodology that 
satisfies the maximum principle and the non-negative 
constraint on general computational grids. 
To this end, we shall use the symbols $\preceq$ 
and $\succeq$ to denote component-wise inequalities 
for vectors. That is, for given any two vectors 
$\boldsymbol{a}$ and $\boldsymbol{b}$ 
%-----------------------------------------;
%  Equation: Definition of preceq symbol  ;
%-----------------------------------------;
\begin{align}
  \boldsymbol{a} \preceq \boldsymbol{b} \quad 
  \mbox{means that } \quad a_i \leq b_i \; \forall i
\end{align}
The symbol $\succeq$ can be similarly defined as well.
Let $<\cdot;\cdot>$ denote the standard inner-product 
in Euclidean space. After finite element discretization, 
the discrete equations corresponding to equation 
\eqref{Eqn:FE_formulation} take the form
%--------------------------------------------------;
%  Equation: Discretization equation for Helmholtz ;
%--------------------------------------------------;
\begin{align}
  \label{Eqn:Helmholtz_discrete}
  \boldsymbol{K} \boldsymbol{c} = \boldsymbol{f}
\end{align}
where $\boldsymbol{K}$ is a symmetric positive definite matrix, $\boldsymbol{c}$ is the 
vector containing nodal concentrations, and $\boldsymbol{f}$ is the force vector. 
Equation~\eqref{Eqn:Helmholtz_discrete} is equivalent to the following minimization problem 
%------------------------------------------------;
%  Equation: Minimization problem for Helmholtz  ;	
%------------------------------------------------;
\begin{align}
  \label{Eqn:Helmholtz_minimization}
  \mathop{\mbox{minimize}}_{\boldsymbol{c} \in \mathbb{R}^{ndofs}} \quad  \frac{1}{2} 
  \langle\boldsymbol{c}; \boldsymbol{K}  \boldsymbol{c}\rangle 
  - \langle\boldsymbol{c}; \boldsymbol{f}\rangle
\end{align}
where ``$ndofs$'' denotes the number of degrees of freedom 
for concentration. Equation~\eqref{Eqn:Helmholtz_discrete} 
can lead to unphysical negative solutions.

Following \cite{Nagarajan_IJNMF_2011_v67_p820,
liska2008enforcing}, a methodology corresponding to equation 
\eqref{Eqn:Helmholtz_minimization} that satisfies the 
non-negative constraint can be written as follows: 
%--------------------------------------;
%  Equation: Non-negative formulation  ;
%--------------------------------------;
\begin{subequations}
  \label{Eqn:non-negative}
  \begin{align}
    &\mathop{\mbox{minimize}}_{\boldsymbol{c} \in \mathbb{R}^{ndofs}} \quad  
    \frac{1}{2} <\boldsymbol{c}; \boldsymbol{K}  \boldsymbol{c}>  - 
    <\boldsymbol{c}; \boldsymbol{f}> \\
    \label{Eqn:non-negative_constraint}
    &\mbox{subject to} \quad \boldsymbol{0} \preceq \boldsymbol{c}  
  \end{align}
\end{subequations}
where $\boldsymbol{0}$ is a vector of size $ndofs$ 
containing zeros. Since $\boldsymbol{K}$ is positive 
definite, equation \eqref{Eqn:non-negative} has a 
unique global minimum \cite{Boyd_convex_optimization}. 
Several robust numerical methods can be used to solve 
equation \eqref{Eqn:non-negative}, which include active 
set strategy, interior point methods 
\cite{Boyd_convex_optimization}. In this paper, to 
solve the resulting optimization problems, we shall 
use the parallel optimization toolkit TAO 
\cite{tao-user-ref}, which has the active-set Newton trust 
region (TRON) and quasi-Newton-based bounded limited memory
variable metric (BLMVM) algorithms.

%==========================================================;
%  Subsection: Governing equations for transient response  ;
%==========================================================;
\subsection{Governing equations for transient response}
We shall denote the time by $t \in [0,\mathcal{I}]$, where 
$\mathcal{I}$ denotes the length of the time interval of 
interest. We shall denote the time-dependent concentration 
by $c(\mathbf{x},t)$. The initial boundary value problem 
can be written as follows:
%-------------------------------------------;
%  Equation: Transient governing equations  ;
%-------------------------------------------;
\begin{subequations}
  \begin{alignat}{2}
    \label{Eqn:LargeScale_GE_unsteady}
    &\frac{\partial c}{\partial t} 
    = \mathrm{div}[\mathbf{D}(\mathbf{x}) \mathrm{grad}[c]] 
    + f(\mathbf{x},t) 
    && \quad \mathrm{in} \; \Omega \times (0,\mathcal{I}) \\
    &c(\mathbf{x},t) = c^{\mathrm{p}}(\mathbf{x},t)  
    && \quad \mathrm{on} \; \Gamma^{\mathrm{D}} 
    \times (0,\mathcal{I}) \\
    -&\widehat{\mathbf{n}}(\mathbf{x}) \cdot 
    \mathbf{D}(\mathbf{x}) \mathrm{grad}[c] 
    = q^{\mathrm{p}}(\mathbf{x},t) 
    && \quad \mathrm{on} \;  \Gamma^{\mathrm{N}} 
    \times (0,\mathcal{I}) \\ 
    \label{Eqn:LargeScale_IC}
    &c(\mathbf{x},0) = c_{0}(\mathbf{x}) 
    && \quad \mathrm{in} \; \Omega 
  \end{alignat}
\end{subequations}
where $c_0(\mathbf{x})$ is the prescribed initial 
concentration, $f(\mathbf{x},t)$ is the time-dependent 
volumetric source/sink, $c^{\mathrm{p}}(\mathbf{x},t)$ is 
the time-dependent prescribed concentration on the 
boundary, and $q^{\mathrm{p}}(\mathbf{x},t)$ is the 
prescribed time-dependent flux on the boundary. 

%====================================================================;
%  Subsubsection: Maximum principle and the non-negative constraint  ;
%====================================================================;
\subsubsection{Maximum principle and the non-negative constraint}
The maximum principle of a transient diffusion equation 
asserts that the maximum can occur only on the boundary 
of the domain or in the initial condition if $f(\mathbf{x},
t) \leq 0$ and $\Gamma^{\mathrm{D}} = \partial \Omega$. 
Mathematically, a solution to equations 
\eqref{Eqn:LargeScale_GE_unsteady}--\eqref{Eqn:LargeScale_GE_unsteady} 
will satisfy:
\begin{align}
  c(\mathbf{x},t) \leq \max \left[\max_{\mathbf{x} \in \Omega} 
    c_0(\mathbf{x}), \max_{\mathbf{x} \in \partial \Omega} 
    c^{\mathrm{p}}(\mathbf{x},t)\right] \quad \forall t
\end{align}
provided $f(\mathbf{x},t) \leq 0$. Similarly, the 
minimum will occur either on the boundary or in the 
initial condition if $f(\mathbf{x},t) \geq 0$. That 
is, if $f(\mathbf{x},t) \geq 0$ then a solution to 
equations \eqref{Eqn:LargeScale_GE_unsteady}--\eqref{Eqn:LargeScale_GE_unsteady} 
satisfies:
\begin{align}
  c(\mathbf{x},t) \geq \min \left[\min_{\mathbf{x} \in \Omega} 
    c_0(\mathbf{x}), \min_{\mathbf{x} \in \partial \Omega} 
    c^{\mathrm{p}}(\mathbf{x},t)\right] \quad \forall t
\end{align}

If $f(\mathbf{x},t) \geq 0$ in $\Omega$, $c^{\mathrm{p}}
(\mathbf{x},t) \geq 0$ on the entire $\partial \Omega$, 
and $c_0(\mathbf{x}) \geq 0$ in $\Omega$ then the maximum 
principle implies that $c(\mathbf{x},t) \geq 0$ in the 
entire domain  at all times, which is the non-negative 
constraint for the concentration field for transient 
problems. 

%=================================================================;
%  Subsubsection: A methodology to enforce the maximum principle  ;
%=================================================================;
\subsubsection{A methodology to enforce the 
  maximum principle for transient problems}
We divide the time interval of interest into 
$\mathcal{N}$ sub-intervals. That is, 
\begin{align}
  [0,\mathcal{I}] := \bigcup_{n = 0}^{\mathcal{N}} [t_n,t_{n+1}]
\end{align}
where $t_n$ denotes the $n$-th time-level. We 
assume that the time-step is uniform, which 
can be written as:
\begin{align}
  \Delta t = t_{n+1} - t_{n}
\end{align}
Following the recommendation provided in  
\cite{Nakshatrala_Nagarajan_Shabouei_Arxiv_2013} to 
meet maximum principles, we employ the backward Euler 
method for temporal discretization. We shall denote 
the nodal concentrations at the $n$-th time-level 
by $\boldsymbol{c}^{(n)}$. We shall denote the minimum 
and maximum values for the concentration by $c_{\mathrm{min}}$ 
and $c_{\mathrm{max}}$, which will be provided by the maximum 
principle and the non-negative constraint. At each time-level, 
one has to solve the following convex quadratic program:
\begin{subequations}
\label{Eq:non-negative}
  \begin{align}
    &\mathop{\mathrm{minimize}}_{\boldsymbol{c}^{(n+1)}} \quad 
    \frac{1}{2} \langle {\boldsymbol{c}^{(n+1)}}; 
    \widetilde{\boldsymbol{K}} \boldsymbol{c}^{(n+1)} \rangle 
    - \langle {\boldsymbol{c}^{(n+1)}}; 
    \widetilde{\boldsymbol{f}}^{(n+1)} \rangle \\
    &\mbox{subject to} \quad c_{\mathrm{min}} \mathbf{1}  \preceq 
    \boldsymbol{c}^{(n+1)} \preceq c_{\mathrm{max}} \mathbf{1}
  \end{align} 
\end{subequations}
where
\begin{align}
  &\widetilde{\boldsymbol{K}} := \frac{1}{\Delta t} \boldsymbol{M} + 
  \boldsymbol{K} \\
  &\widetilde{\boldsymbol{f}}^{(n+1)} := \boldsymbol{f}^{(n+1)} + 
  \frac{1}{\Delta t} \boldsymbol{M} \boldsymbol{c}^{(n+1)}
\end{align}
In the above equation, $\boldsymbol{M}$ is the capacity 
matrix \cite{Nakshatrala_Nagarajan_Shabouei_Arxiv_2013}.

%*********************************************;
%                                             ;
%  NAME                                       ;
%    S3_LargeScale_PI.tex          	      ;
%                                             ;
%  WRITTEN BY                                 ;
%    Justin Chang			      ;
%    Satish Karra			      ;
%    Kalyana Babu Nakshatrala                 ;
%                                             ;
%*********************************************;
\section{PARALLEL IMPLEMENTATION}
\label{Sec:Parallel_Implementation}
\subsection{PETSc and TAO}
%------------------------------------------------;
%  Figure: DMPlex data structure                 ;
%------------------------------------------------;
\begin{figure}[htp]
\centering
\subfloat[Optimal 2D and 3D elements]{\includegraphics[scale=0.75]{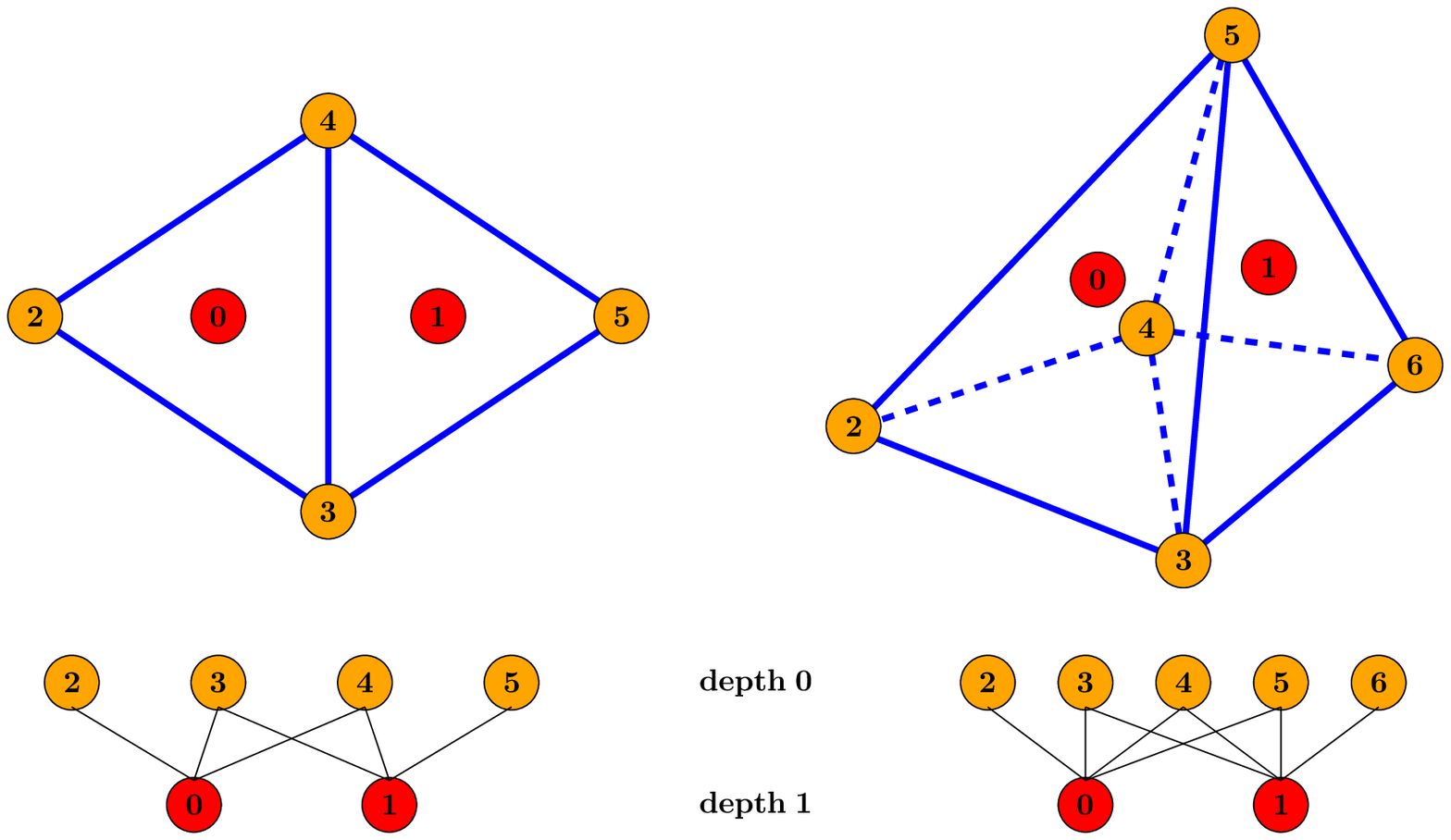}}
\label{Fig:S3_dmplex_simple}\\
\subfloat[Interpolated 2D and 3D elements]{\includegraphics[scale=0.75]{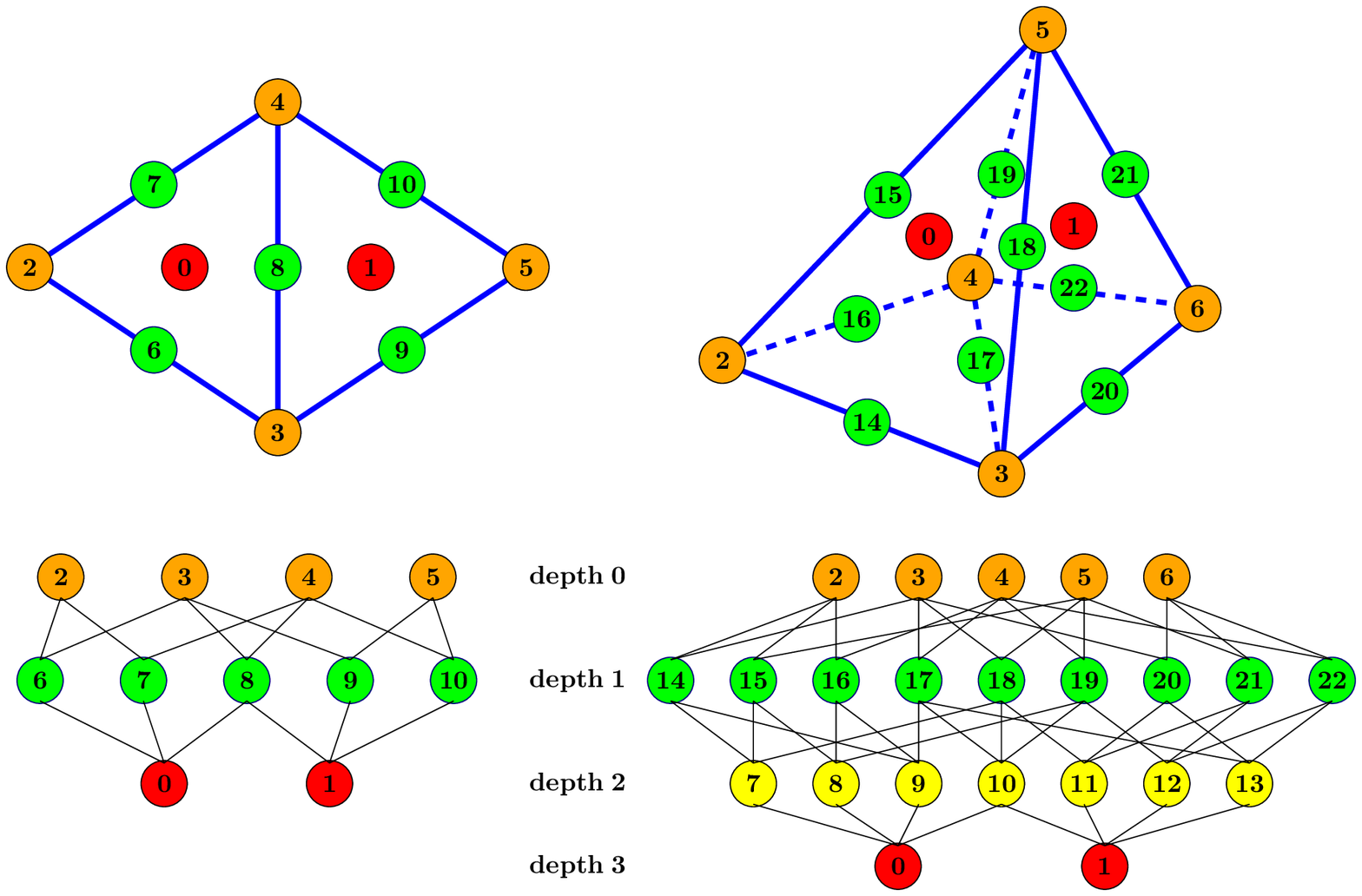}}
\label{Fig:S3_dmplex_interpolated}
\caption{Representation of mesh points within the DMPlex data structure and their associated directed acyclic graphs}
\label{Fig:S3_dmplex}
\end{figure}
We leverage on existing scientific libraries such as PETSc and TAO to formulate our 
large-scale computational framework. PETSc is a suite of data structures and 
routines for the parallel solution of scientific applications. It also provides 
interfaces to several other libraries such as Metis/ParMETIS\cite{METIS} and HDF5 
\cite{hdf5} for mesh partitioning and binary data format handling respectively. The 
Data Management (DM) data structure is used to manage all information including 
vectors and sparse matrices and compatible with binary data formats. To handle 
unstructured grids in parallel, a subset of the DM structure called DMPlex 
(see \cite{Lange_easc,Knepley_scicomp,petsc-user-ref}), as shown in Figure \ref{Fig:S3_dmplex},
uses the direct acyclic graph to organize  all mesh information. This topology 
enables the freedom to mix and match various non vertex-based discretization such as
the two-point flux finite volume method and the classical mixed formulations
based on the lowest-order Raviart Thomas finite element space.

Another important feature within PETSc is TAO. The TAO library has a suite of data 
structures and routines that enable the solution of large-scale optimization 
problems. It can support any data structure or solver within PETSc. Our non-negative
methodology will use both the Newton-Trust Region (TRON) and Bounded Limited-Memory 
Variable-Metric (BLMVM) solvers available within TAO. BLMVM is a quasi-Newton method 
that uses projected gradients to approximate the Hessian, which is useful for problems where the
hessian is too complicated or expensive to compute. Other optimization algorithms 
such as TRON and the Gradient Projected Conjugate Gradient (GPCG) typically require 
Hessian information and more memory, but they are expected to converge more rapidly
than BLMVM. Further details regarding the implementation of these various methods 
may be found in \cite{tao-user-ref} and the references within.
 
%======================================================;
%  Subsection: Implementation ;
%======================================================;
\subsection{Finite element implementation}
PETSc abstractions for finite elements, quadrature rules, and function spaces have 
also been recently introduced and are suitable for the mesh topology within DMPlex. 
They are built upon the same framework as the Finite element Automatic Tabulator 
(FIAT) found within the FEniCS Project \cite{Kirby2012a,Alogg,LoggMardalEtAl2012a}. 
The finite element discretizations simply need the equations, auxiliary coefficients
(e.g., permeability, diffusivity, etc.), and boundary conditions specified as 
point-wise functions. We express all discretizations in nonlinear form so let 
$\boldsymbol{r}$ and $\boldsymbol{J}$ denote the residual and Jacobian respectively.

Following the FEM model outlined in \cite{2013arXiv1309.1204K}, we consider the weak
form that depends on fields and gradients. The residual evaluation can be expressed
as: 
\begin{align}
&\boldsymbol{w}^{\mathrm{T}}\boldsymbol{r}(\boldsymbol{c}) \sim \int_{\Omega^{e}} \left[w\cdot 
\mathcal{F}_{0}\left(c,\nabla c\right) + 
\nabla w \cdot \boldsymbol{\mathcal{F}}_1(c,\nabla c) 
\right]\mathrm{d} \Omega = 0
\end{align}
where $\mathcal{F}_0(c,\nabla c)$ and $\boldsymbol{\mathcal{F}}_1(c,\nabla c)$ are 
user-defined point-wise functions that capture the problem physics. This framework 
decouples the problem specification from the mesh and degree of freedom traversal. 
That is, the scientist need only focus on providing point function evaluations while
letting the finite element library take care of meshing, quadrature points, basis 
function evaluation, and mixed forms if any. The discretization of the residual is 
written as: 
\begin{align}
&\boldsymbol{r}(\boldsymbol{c}) = \mathop{\huge \boldsymbol{\mathsf{A}}}_{e=1}^{Nele} 
\left[\begin{array}{lr}\boldsymbol{N}^{\mathrm{T}} & \boldsymbol{B}^{\mathrm{T}}	
\end{array}\right]
\boldsymbol{W}\left[\begin{array}{l}	\mathcal{F}_0(c_q,\nabla c_q) \\ 
\boldsymbol{\mathcal{F}}_1(c_q,\nabla c_q)\end{array}\right]
\label{Eqn:S3_residual}
\end{align}
where $\boldsymbol{\mathsf{A}}$ represents the standard assembly operator, 
$\boldsymbol{N}$ and $\boldsymbol{B}$ are matrix forms of basis functions that 
reduce over quadrature points, $\boldsymbol{W}$ is a diagonal matrix of quadrature 
weights (including the geometric Jacobian determinant of the element), and $c_q$ is 
the field value at quadrature point $q$. The Jacobian of \eqref{Eqn:S3_residual} 
needs only the derivatives of the point-wise functions:
\begin{align}
&\boldsymbol{J}(\boldsymbol{c}) = \mathop{\huge \boldsymbol{\mathsf{A}}}_{e=1}^{Nele} 
\left[\begin{array}{lr}\boldsymbol{N}^{\mathrm{T}} & \boldsymbol{B}^{\mathrm{T}}	\end{array}\right]
\boldsymbol{W}\left[\begin{array}{ll} \mathcal{F}_{0,0} & \mathcal{F}_{0,1} \\ 
\boldsymbol{\mathcal{F}}_{1,0}& \boldsymbol{\mathcal{F}}_{1,1}\end{array}\right]
\left[\begin{array}{c} \boldsymbol{N} \\ \boldsymbol{B}
\end{array}\right], \quad \left[\mathcal{F}_{i,j}\right] = 
\left[\begin{array}{ll}
\frac{\partial\mathcal{F}_0}{\partial c} & \frac{\partial\mathcal{F}_0}{\partial\nabla c} \\
\frac{\partial\boldsymbol{\mathcal{F}}_1}{\partial c} & \frac{\partial\boldsymbol{\mathcal{F}}_1}{\partial\nabla c}
\end{array}\right](c_q,\nabla c_q)
\label{Eqn:S3_jacobian}
\end{align}
The point-wise functions corresponding to the weak form in 
\eqref{Eqn:functionals_B_L} would be:
%-------------------------;
%  Algorithm: Pseudocode  ;
%-------------------------;
\begin{algorithm}[b]
  \begin{algorithmic}
    \State Create/input DAG on rank 0
    \State Create/input cell-wise velocity on rank 0
    \If {size $>$ 1}
    \State Partition mesh among all processors
    \EndIf
    \State Refine distributed mesh if necessary
    \State Create PetscSection and FE discretization
    \State Set $n=0$ and $\boldsymbol{c}^{(0)}=10^{-8} $
    \State Insert Dirichlet BC constraints into  $\boldsymbol{c}^{(0)} $
    \State Compute Jacobian $\boldsymbol{J}$
    \While {true} \Comment{Begin time-stepping scheme}
	    \State Compute Residual $\boldsymbol{r}^{(n)}$ 
	    \If {Classical Galerkin}
	    \Comment{Solve without non-negative methodology}
	    \State $\boldsymbol{c}^{(n+1)} = \boldsymbol{c}^{(n)} - \boldsymbol{J}\backslash\boldsymbol{r}^{(n)}$
	    \Else\Comment{Solve with non-negative methodology}
	    \State TaoSolve() for $\boldsymbol{c}^{(n+1)}$ based on equations \eqref{Eqn:S3_tao_objective} and \eqref{Eqn:S3_tao_gradient}
	    \EndIf
	    \If {steady-state or $(n) == $ total number of time steps}
	    \State break
	    \Else
	    \State $n += 1$
	    \EndIf
    \EndWhile
    
  \end{algorithmic}
  \caption{Pseudocode for the large-scale transport solver \label{Algo:S3_pseudocode}}
\end{algorithm}
\begin{subequations}
\begin{align}
&\mathcal{F}_0 = -f(\mathbf{x}),\quad \boldsymbol{\mathcal{F}}_1 = \mathbf{D}(\mathbf{x})
\nabla c_q \\
&\mathcal{F}_{0,0} = 0, \quad \mathcal{F}_{0,1} = 0, \quad \boldsymbol{\mathcal{F}}_{1,0} = \mathbf{0}, \quad \boldsymbol{\mathcal{F}}_{1,1} = \mathbf{D}(\mathbf{x})
\end{align}
\end{subequations}
Similarly, the point-wise functions for the transient response are:
\begin{subequations}
\begin{align}
&\mathcal{F}_0 = \dot{c}_q-f(\mathbf{x},t),\quad \boldsymbol{\mathcal{F}}_1 = \mathbf{D}(\mathbf{x})
\nabla c_q \\
&\mathcal{F}_{0,0} = \frac{1}{\Delta t}, \quad \mathcal{F}_{0,1} = 0, \quad \boldsymbol{\mathcal{F}}_{1,0} = \mathbf{0}, \quad \boldsymbol{\mathcal{F}}_{1,1} = \mathbf{D}(\mathbf{x})
\end{align}
\end{subequations}
where $\dot{c}_{q}$ denotes the time derivative. A similar discretization is used 
to project the Neumann boundary conditions into the residual vector. 
Assuming a fixed time-step, $[\mathcal{F}_{i,j}]$ and the Jacobian in equation 
\eqref{Eqn:S3_jacobian} do not change with time and have to be computed only once. 
If $n$ denotes the time level ($n$ = 0 denotes initial condition) then 
the residual and Jacobian can be defined as:
\begin{align}
&\boldsymbol{r}^{(n)} \equiv \boldsymbol{r}(\boldsymbol{c}^{(n)})\\
&\boldsymbol{J} \equiv \boldsymbol{J}(\boldsymbol{c}^{(0)})
\end{align}
To enforce the non-negative methodology, the following objective function $b$ 
and gradient function $\boldsymbol{g}$ is provided: 

%------------------------------------------------;
%  Equation: TAO routines                        ;
%------------------------------------------------;
  \begin{align}
  \label{Eqn:S3_tao_objective}
    &b = \frac{1}{2}\boldsymbol{c}^{(n+1)}\cdot\boldsymbol{J}\boldsymbol{c}^{(n+1)} + \boldsymbol{c}^{(n+1)}\cdot\left[\boldsymbol{r}^{(n)}-\boldsymbol{J}\boldsymbol{c}^{(n)}\right]\\
    \label{Eqn:S3_tao_gradient}
    &\boldsymbol{g} = \boldsymbol{J}\left[\boldsymbol{c}^{(n+1)}-\boldsymbol{c}^{(n)}\right] + \boldsymbol{r}^{(n)}
  \end{align}
BLMVM relies only on the above two equations, whereas TRON needs the Hessian which 
is equivalent to $\boldsymbol{J}$. Algorithm \ref{Algo:S3_pseudocode} outlines 
the steps taken in our computational framework.

%*********************************************;
%                                             ;
%  NAME                                       ;
%    S4_LargeScale_PM.tex          	      ;
%                                             ;
%  WRITTEN BY                                 ;
%    Justin Chang			      ;
%    Satish Karra			      ;
%    Kalyana Babu Nakshatrala                 ;
%                                             ;
%*********************************************;
\section{PERFORMANCE MODELING}
\label{Sec:Performance_Modeling}
PETSc is a constantly evolving open-source library that
brings out new features and algorithms almost every day.
It has capabilities to interface with a large number of
other open-source software and linear algebra packages. 
However, it is not always known which of these algorithms 
will have the best performance across multiple distributed 
memory HPC systems, especially if these packages have little 
documentation and have to be used as black-box solvers. 
Computational scientists would like to know which solvers 
or algorithms to use for their specific need before running 
jobs on the state-of-the-art HPC systems. The first and the 
trivial metric to look for in answering this question is the 
time-to-solution for a given solver or optimization method. 
However, additional information is needed in order to quantify 
the hardware and algorithmic efficiency as well as the potential 
scalability across multiple cores in the strong sense.
%------------------------------------------------;
%  Table: HPC specification                      ;
%------------------------------------------------;
\begin{table}[htp]
  \centering
    \caption{List of HPC systems used in this study \label{Tab:S3_HPC}}
  \begin{tabular}{lcc}
    \hline
    & Mustang (MU) & Wolf (WF) \\
    \hline
    Processor & AMD Opteron 6176 & Intel Xeon E5-2670 \\
    Clock rate& 2.3 GHz & 2.6 GHz\\
    FLOPs/cycle & 4 & 8\\
    Sockets per compute node & 2 & 2\\
    NUMA nodes per socket & 2 & 1\\
    Cores per socket & 12 & 8\\
    Total cores (compute nodes) & 38400 (1600) & 9856 (616)\\
    Memory per compute node & 64 GB & 64 GB\\
    L1 cache per core & 128 KB &32 KB\\
    L2 cache per core & 512 KB & 256 KB\\
    L3 cache per socket & 12 MB & 20 MB\\
    Interconnect & 40 Gb/s & 40 Gb/s \\
    \hline
  \end{tabular}
\end{table}

Hardware specifications of HPC systems significantly impact the performance of any 
numerical algorithm. Ideally we want our simulations to consume as little wall-clock
time as possible as the number of processing cores increases (i.e., achieving good 
speedup), but several other factors including compiler vectorization, cache locality, 
memory bandwidth, and code implementation may drastically affect the performance. 
Table \ref{Tab:S3_HPC} lists the hardware specifications 
of the two HPC systems (Mustang and Wolf) that are used 
in our numerical experiments. The Mustang HPC system consists of
relatively older generation of processors so  
it is expected to not perform as well. One could 
simply measure wall-clock time across multiple compute 
nodes on the respective HPC machines and determine the 
parallel efficiency of a certain algorithm, but we are 
interested in quantifying how different algorithms behave 
sequentially and what kind of parallel performance to expect 
before running numerical simulations on supercomputers. The wall-clock 
time of any simulation can generally be summed up as a function of three things: 
the workload, transfer of data between the memory and CPU register, 
and interprocess communication. Hardware efficiency in this 
context is defined as the amount of time spent performing work 
over waiting on memory fetching and cache registers to free up. 

The limiting factor of performance for numerical methods on modern computing 
architectures is upper-bounded by the memory bandwidth. That is, the floating point 
performance given by FLOPS/s will never reach the theoretical peak performance 
(TPP). This limitation is particularly important for iterative solvers and 
optimization methods that rely on numerous sparse matrix-vector (SpMV) 
multiplications (see \cite{may_ptatin} and the references within). The 
frequent use of SpMV allows for little cache reuse and will result in a 
large number of very expensive cache misses. Such behavior is important to 
document when determining how efficient a scientific code is, so performance models 
such as the Roofline Model \cite{williams,Lo_roofline}, which measures memory transfers, 
have been used to better quantify the efficiency with respect to the hardware. 
Performance models in general can help application developers 
identify bottlenecks and indicate which areas of the code can be further optimized. 
In other words, the code can be designed so that it maximizes the full 
benefits of the available computing resources. In the next 
section, we will demonstrate that such models can also be used to predict the 
parallel efficiency of various optimization solvers on the two very different LANL 
HPC systems. The key parameter for these performance models is the Arithmetic 
Intensity (AI) which is defined as:
%------------------------------------------------;
%  Equation: Arithmetic Intensity                ;
%------------------------------------------------;
\begin{align}
   &\mbox{AI} = \frac{\mbox{Total FLOPS}}
   {\mbox{Total Bytes Transferred}}
  \label{Eqn:S3_efficiency}
\end{align}
where the Total Bytes Transferred (TBT) metric denotes
the amount of bandwidth needed for a given floating
point operation. The AI serves as a multiplier to the
actual memory bandwidth and creates a ``roofline" for
the estimation of ideal peak performance. A cache model 
is needed in order to properly define the TBT. 

%------------------------------------------------;
%  Table: AI chart                               ;
%------------------------------------------------;
\begin{table}[tp]
  \centering
    \caption{Commonly used PETSc operations and their 
      respective Total Bytes Transferred. Here we note 
      $X,Y,Z$ as vectors with $i = 1,\cdots,N$ entries, 
    $a$ is a scalar value, and $nz$ denotes the total 
      number of non-zeros. We assume that the sizes of 
      integers and doubles are 4 and 8 bytes respectively. 
    \label{Tab:S3_TBT}}
  \begin{tabular}{lccc}
    \hline
    PETSc function & Operation & Total Bytes Transferred \\
    \hline
    VecNorm() & $a = \sqrt{\sum^{N}_{i} X(i)^2} $ & $8(N + 1)$ \\
    VecDot() & $a = \sum^{N}_{i}X(i)*Y(i)$ & $8(2N + 1)$ \\
    VecCopy() & $Y \leftarrow X$ & $8(2N)$ \\
    VecSet() & $Y(i) = a$ & $8(2N)$ \\
    VecScale() & $Y = a*Y$ & $8(2N)$ \\
    VecAXPY() & $Y = a*X+Y$ & $8(3N)$ \\
    VecAYPX() & $Y = X+a*Y$ & $8(3N)$ \\
    VecPointwiseMult() & $Z(i) = X(i)*Y(i)$ & $8(3N)$ \\
    MatMult() & SpMV & $4(N + nz) + 8(2N + nz)$ \\
    \hline
  \end{tabular}
\end{table}
To this end, we propose a roofline-like performance 
model where the TBT assumes a ``perfect cache" -- each 
byte of the data needs to be fetched from DRAM only 
once. This assumption enables us to predict a slightly 
more realistic upper bound of the peak performance than 
by simply comparing to the TPP. Table \ref{Tab:S3_TBT} 
lists the key PETSc functions used for the solvers and their 
respective estimates of TBT based on the perfect cache assumption. 
The formula for SpMV follows the procedure outlined in 
\cite{Kaushik99towardrealistic}. We assume that the TBT formula 
for operations also involving a sparse matrix and vector 
like the incomplete lower-upper (ILU) factorization to be 
the same as MatMult(). Estimating the TBT for other 
important operations like the sparse matrix-matrix and 
triple matrix products (which are important for multi-grid methods) 
is an area of future work. In short, our AI formulation 
relies on the following four key assumptions:
\begin{enumerate}[(i)]
\item All floating-point operations (add, multiply, square roots, etc.) are treated 
equally and equate to one FLOP count. 
\item There are no conflict misses. That is, each matrix and vector element is 
loaded into cache only once.
\item Processor never waits on a memory reference. That is, any number of loads and 
stores are satisfied in a single cycle.
\item Compilers are capable of storing scalar multipliers in the register only for 
pure streaming computations.
\end{enumerate}
Therefore, the efficiency based on this new roofline-like performance model as:
%------------------------------------------------;
%  Figure: STREAMS Benchmark                     ;
%------------------------------------------------;
\begin{figure}[htp]
\centering
\subfloat{\includegraphics[scale=0.55]{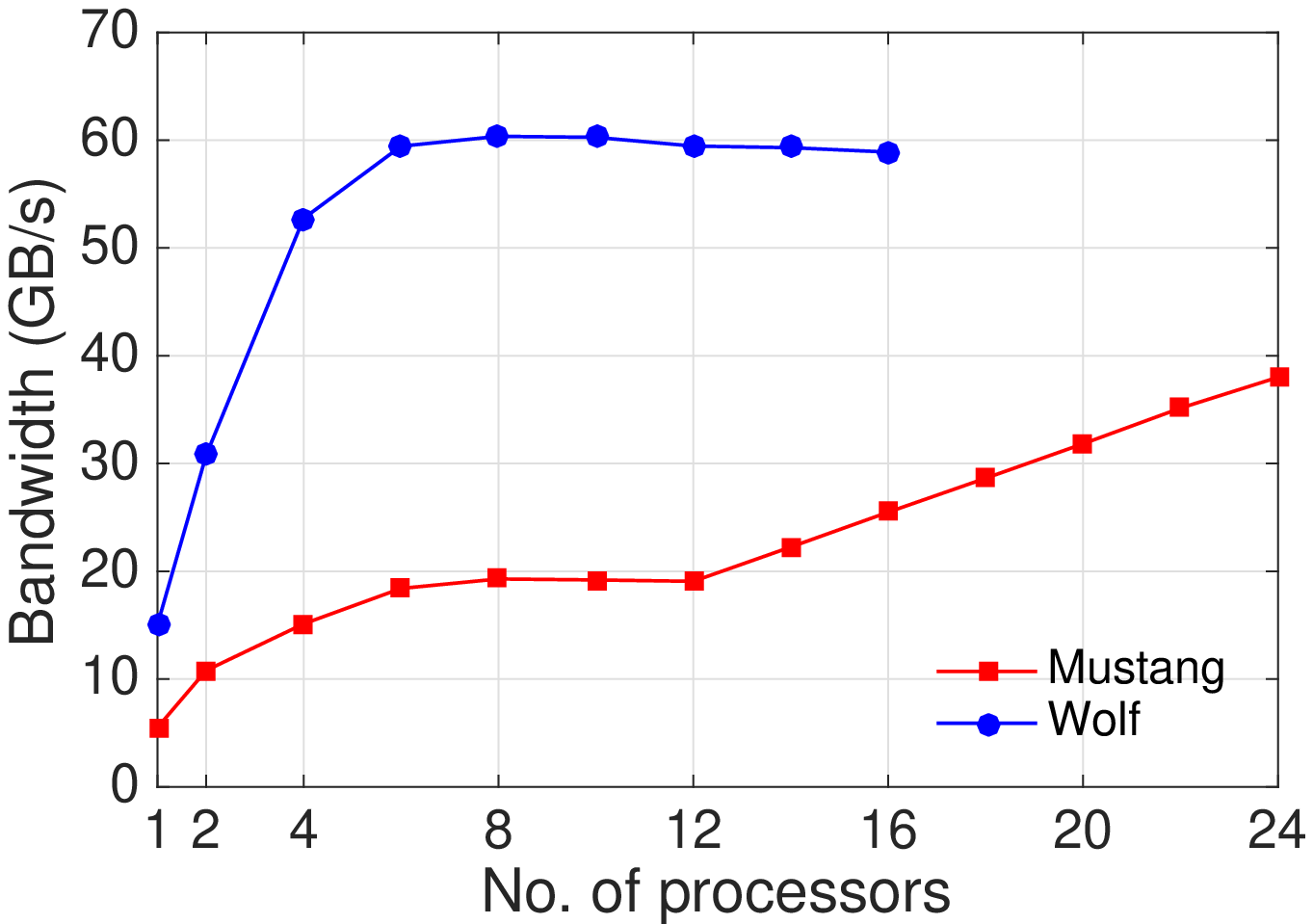}}
\caption{Estimated memory bandwidth of a single Mustang and Wolf compute node based 
on the STREAMS Triad Benchmark}
\label{Fig:S3_streams}
\end{figure}
%%------------------------------------------------;
%%  Equation: Ideal flopss                        ;
%%------------------------------------------------;
\begin{align}
   &\mbox{Efficiency (\%)} = \frac{\mbox{Measured FLOPS/s}}
   {\mathrm{min}\left\{\begin{array}{l}
    \mbox{TPP} \\
    \mathrm{AI}\times\mbox{STREAMS}
  \end{array}
  \right.}\times100
  \label{Eqn:S3_efficiency}
\end{align}
where the numerator is reported by the PETSc program
and the denominator is the ideal performance upper-bounded 
by both the TPP and the product of AI and STREAMS bandwidth. 
STREAMS Triad \cite{streams} is one of the most 
popular benchmarks for determining the achievable memory 
performance of a given machine. Figure \ref{Fig:S3_streams} denotes the 
estimated  memorybandwidth as a function of number of cores on a single 
Mustang and Wolf node. It is interesting to note that although
the Wolf node has a greater bandwidth, there is
no performance gain past eight cores. This means that an optimal use of a
Wolf compute node for memory-bandwidth bound algorithms would be 
eight cores, whereas one would still see some performance gains when
using all 24 cores on a Mustang node.

The performance model that uses equation \eqref{Eqn:S3_efficiency} 
is a serial model so the STREAMS metric for the Mustang and Wolf 
systems are 5.65 GB/s and 15.5 GB/s respectively. It should be 
noted that this performance model does account for cache effects. 
That is, it does not quantify the useful bandwidth sustained for 
some level of cache. The true hardware and algorithmic efficiency 
is not be reflected by this model, so our aim is to show relative 
performance between select PETSc and TAO solvers. Comparing the 
AI and the measured FLOPS/s with the STREAMS bandwidth will give 
us a better understanding of how high-performing the PETSc and TAO 
solvers are for select problems.

%*********************************************;
%                                             ;
%  NAME                                       ;
%  S5_LargeScale_NR.tex          	      ;
%                                             ;
%  WRITTEN BY                                 ;
%	 Justin Chang		              ;
%	 Satish Karra			      ;
%  Kalyana Babu Nakshatrala                   ;
%                                             ;
%*********************************************;
%------------------------------------------------;
%  Figure: Unit cube description                 ;
%------------------------------------------------;
\begin{figure}[t]
\centering
\subfloat[Location of the hole]{\includegraphics[scale=0.35]{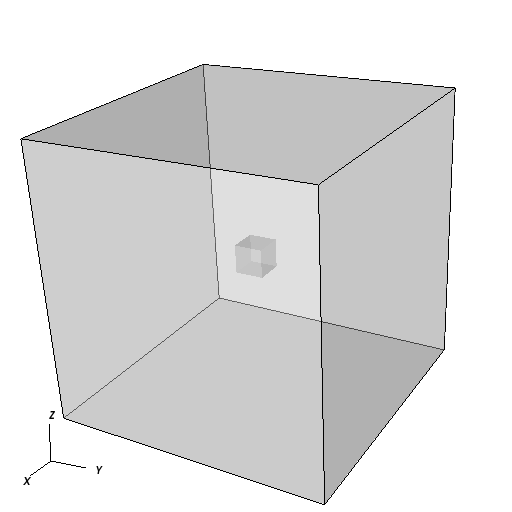}}
\subfloat[Mesh type A]{\includegraphics[scale=0.35]{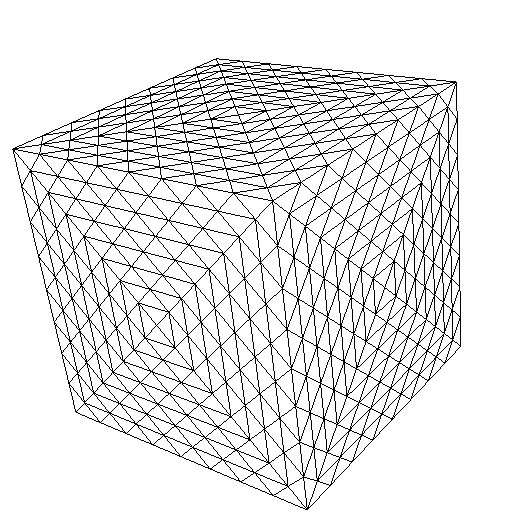}}\\
\subfloat[Mesh type B]{\includegraphics[scale=0.35]{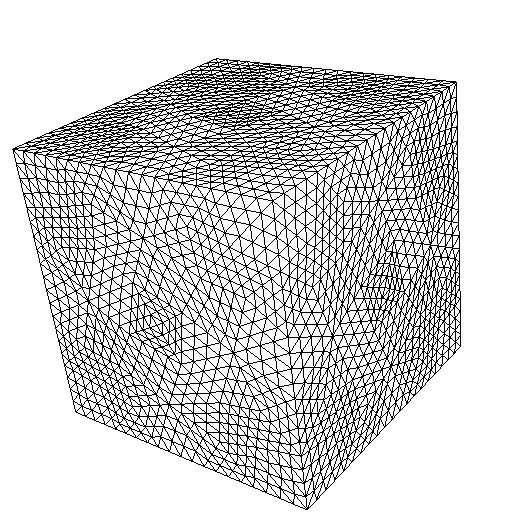}}
\subfloat[Mesh type C]{\includegraphics[scale=0.35]{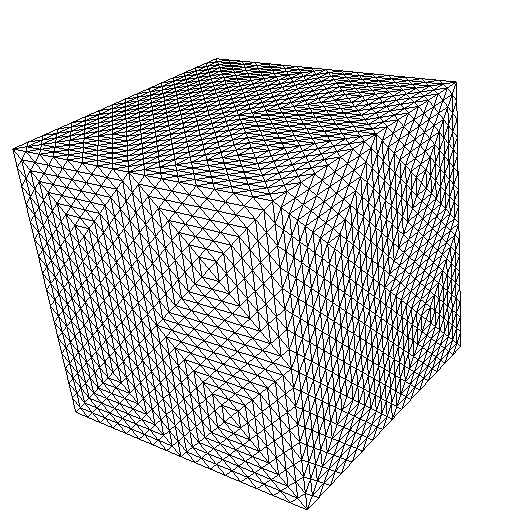}}
\caption{Cube with a hole: pictorial description and the associated unstructured grids.}
\label{Fig:S5_cube_description}
\end{figure}
%------------------------------------------------;
%  Figure: Numerical results of unit cube        ;
%------------------------------------------------;
\begin{figure}[t]
\centering
\subfloat{\includegraphics[scale=0.55]{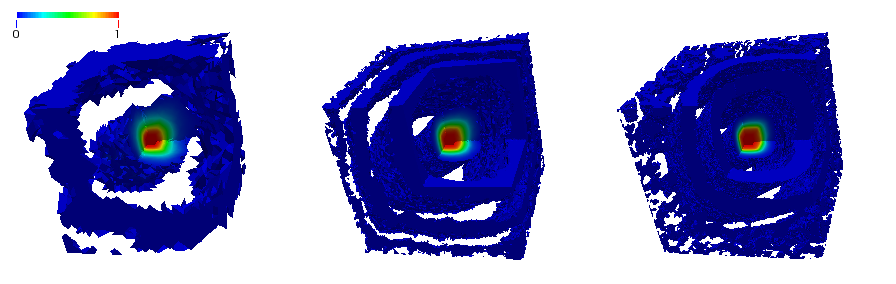}}\\
\subfloat{\includegraphics[scale=0.55]
{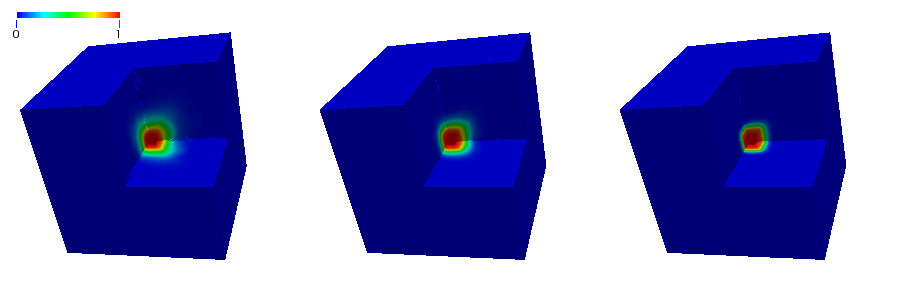}}
\caption{Cube with a hole: numerical solution for cases A1 (left), 
B2 (middle), and C3 (right) using the Galerkin formulation (top 
row) and non-negative methodology (bottom row).}
\label{Fig:S5_cube_results}
\end{figure}
\section{REPRESENTATIVE NUMERICAL RESULTS}
\label{Sec:Numerical_Results}
In this section, we compare the performance of our 
non-negative methodology using the TAO solver to 
that of the Galerkin formulation using the Krylov 
Subspace (KSP) solver. We examine the performance 
using two problems: 
\begin{enumerate}[(i)]
\item a unit cube with a hole under steady-state, and 
\item a transient Chromium transport problem.
\end{enumerate}
The diffusivity tensor is assumed to depend 
on the flow velocity through 
\begin{align}
  \mathbf{D}(\mathbf{x}) = \left(\alpha_T \|\mathbf{v}\|+D_{M}\right) \mathbf{I} + 
(\alpha_L - \alpha_T) \frac{\mathbf{v} \otimes \mathbf{v}}{\|\mathbf{v}\|} 
\label{Eqn:S5_dispersion_tensor}
\end{align}
where $\alpha_L$, $\alpha_T$, and $D_M$ denote the longitudinal dispersivity,
transverse dispersivity and molecular diffusivity, respectively. We employ the 
conjugate gradient method and the block Jacobi/ILU(0) preconditioner for solving 
the linear system from the Galerkin formulation and employ TAO's TRON and BLMVM 
methods for the non-negative methodology. The relative convergence tolerances for 
both KSP and TAO solvers are set to $10^{-6}$, and $\Delta t$ for the transient 
response in the Chromium problem is initially set to 0.2 days. For strong-scaling 
studies shown here, we used OpenMPI v1.6.5 for message passing and bound processes to cores 
while mapping by sockets. ParaView\cite{paraview} and VisIt\cite{VisIt} were 
used to generate all contour and mesh plots.

\begin{remark}
Throughout the paper, the non-negative methodology that we refer to, is in fact a 
discrete maximum principle preserving methodology, in that, along with the 
non-negative constraint we also enforce that the concentrations are less than 
or equal to 1.
\end{remark}

%====================================================;
%  Subsection: Anisotropic diffusion in a unit cube  ;
%====================================================;
\subsection{Anisotropic diffusion in a unit cube with a cubic hole}
Let the computational domain be a unit cube with a cubic hole of size 
$[4/9, 5/9] \times [4/9, 5/9] \times [4/9, 5/9]$. The concentration on the 
outer boundary is taken to be zero and the concentration on the interior 
boundary is taken to be unity. The volumetric source is taken as zero 
(i.e., $f(\mathbf{x}) = 0$). The velocity vector field for this problem is 
chosen to be
\begin{align}
\mathbf{v}(\mathbf{x}) = \mathbf{e}_x + \mathbf{e}_y + \mathbf{e}_z
\end{align}
%------------------------------------------------;
%  Table: Unit cube information                  ;
%------------------------------------------------;
\begin{table}[t]
  \centering
  \caption{Cube with a hole: list of various mesh type and refinement 
    level combinations used \label{Tab:S5_cube_information}}
  \begin{tabular}{lcccc}
    \hline
    Case & Mesh type & Refinement level & Tetrahedrons & Vertices \\
    \hline
    A1 & A & 1 & 199,296 & 36,378 \\
    B1 & B & 1 & 409,848 & 75,427 \\
    C1 & C & 1 & 793,824 & 140,190 \\
    A2 & A & 2 & 1,594,368 & 278,194 \\
    B2 & B & 2 & 3,278,784 & 574,524 \\
    C2 & C & 2 & 6,350,592 & 1,089,562 \\
    A3 & A & 3 & 12,754,994 & 2,175,330 \\
    B3 & B & 3 & 26,230,272 & 4,483,126 \\
    C3 & C & 3 & 50,804,736 & 9,172,044 \\
    \hline
  \end{tabular}
\end{table}
%------------------------------------------------;
%  Table: Unit cube DMP                         ;
%------------------------------------------------;
\begin{table}[t]
  \centering
  \caption{Cube with a hole: minimum and maximum concentrations 
    for each case \label{Tab:S5_cube_dmp}}
  \begin{tabular}{lccc}
    \hline
    Case & Min. concentration & Max. concentration & \% nodes violated \\
    \hline
    A1 & -0.0224825 & 1.00000 & 9,518/36,378 $\rightarrow$ 26.2\% \\
    B1 & -0.0139559 & 1.00000 & 32,247/43,180 $\rightarrow$ 42.8\% \\
    C1 & -0.0125979 & 1.00000 & 57,272/140,190 $\rightarrow$ 40.9\% \\
    A2 & -0.0311518 & 1.00103 & 82,983/278,194 $\rightarrow$ 29.2\% \\
    B2 & -0.0143857 & 1.00000 & 255,640/574,524 $\rightarrow$ 44.9\% \\
    C2 & -0.0119539 & 1.00972 & 453,766/1,089,562 $\rightarrow$ 41.6\% \\
    A3 & -0.0258559 & 1.00646 & 643,083/2,175,330 $\rightarrow$ 29.6\% \\
    B3 & -0.0115908 & 1.00192 & 2,073,934/4,483126 $\rightarrow$ 46.3\%  \\
    C3 & -0.0096186 & 1.00545 & 4,932,551/9,172,044 $\rightarrow$ 53.8\%  \\
    \hline
  \end{tabular}
\end{table}
%------------------------------------------------;
%  Table: Unit cube timings Mustang              ;
%------------------------------------------------;
\begin{table}[b]
  \centering
  \caption{Cube with a hole: wall-clock times (seconds) on Mustang for each solver \label{Tab:S5_cube_times_mustang}}
  \begin{tabular}{lccccc}
    \hline
    Case & Galerkin & TRON1 & TRON2 & TRON3 & BLMVM \\
    \hline
    A1 & 0.337 & 0.933 & 0.981 & 1.14 & 2.62 \\
    B1 & 0.790 & 1.72 & 2.06 & 2.71 & 5.04 \\
    C1 & 2.24 & 4.34 & 5.80  & 7.74 & 13.5 \\
    A2 & 7.21 & 15.2 & 21.7 & 32.5 & 72.0 \\
    B2 & 15.4 & 30.0 & 43.7 & 57.5 & 109 \\
    C2 & 40.4 & 67.8 & 113 & 118 & 286 \\
    A3 & 121 & 225 & 414 & 599 & 1167 \\
    B3 & 315 & 498 & 1061 & 1344 & 2524 \\
    C3 & 997 & 1539 & 2490 & 4365 & 9679 \\
    \hline
  \end{tabular}
\end{table}
%------------------------------------------------;
%  Table: Unit cube timings Wolf                 ;
%------------------------------------------------;
\begin{table}[b]
  \centering
  \caption{Cube with a hole: wall-clock times (seconds) on Wolf for each solver \label{Tab:S5_cube_times_wolf}}
  \begin{tabular}{lccccc}
    \hline
    Case & Galerkin & TRON1 & TRON2 & TRON3 & BLMVM \\
    \hline
    A1 & 0.126 & 0.388 & 0.396 & 0.449 & 1.01 \\
    B1 & 0.314 & 0.720 & 0.853 & 1.07 & 2.03 \\
    C1 & 0.888 & 1.91 & 2.47 & 3.31 & 5.71 \\
    A2 & 2.58 & 6.34 & 8.74 & 12.8 & 26.2 \\
    B2 & 5.90 & 12.9 & 17.8 & 22.8 & 46 \\
    C2 & 16.2 & 30.1 & 47.3 & 48.9 & 133 \\
    A3 & 48.0 & 98.4 & 129 & 247 & 609 \\
    B3 & 107 & 171 & 342 & 435 & 1060 \\
    C3 & 281 & 467 & 870 & 1245 & 3131 \\
    \hline
  \end{tabular}
\end{table}

The diffusion parameters are set as: $\alpha_L = 1$, $\alpha_T = 0.001$, and 
$D_M = 0$. The pictorial description of the computational domain and the three mesh
types composed of 4-node tetrahedrons are shown in Figure 
\ref{Fig:S5_cube_description}. We consider three unstructured mesh types 
with three levels of element-wise mesh refinement, giving us nine total case 
studies of increasing problem size as shown in Table \ref{Tab:S5_cube_information}. 
We ran a total of five different simulations for this study: 
\begin{itemize}
\item Galerkin with CG/block Jacobi
\item TRON1: with KSP tolerance of $10^{-1}$
\item TRON2: with KSP tolerance of $10^{-2}$
\item TRON3: with KSP tolerance of $10^{-3}$
\item BLMVM
\end{itemize} 
The TRON solvers also use the CG and block Jacobi preconditioner but with 
different KSP tolerances. Numerical results for both the Galerkin formulation and 
the non-negative methodologies for some of the mesh cases are shown in Figure 
\ref{Fig:S5_cube_results}.  The top row of figures arise from the Galerkin 
formulation where the white regions denote negative concentrations, and the 
bottom row arise from either TRON or BLMVM. Details concerning the 
violation of the DMP for each case study can be found in Table 
\ref{Tab:S5_cube_dmp}. Concentrations both negative
and greater than one arise for all case studies.  Moreover, 
simply refining the mesh does not resolve these issues; in fact, 
refinment worsens the violation. These numerical results indicate that our 
computational framework can successfully enforce the DMP for diffusion problems 
with highly anisotropic diffusivity.
%------------------------------------------------;
%  Figure: Solver iterations - Galerkin          ;
%------------------------------------------------;
\begin{figure}[t]
\centering
\subfloat{\includegraphics[scale=0.55]{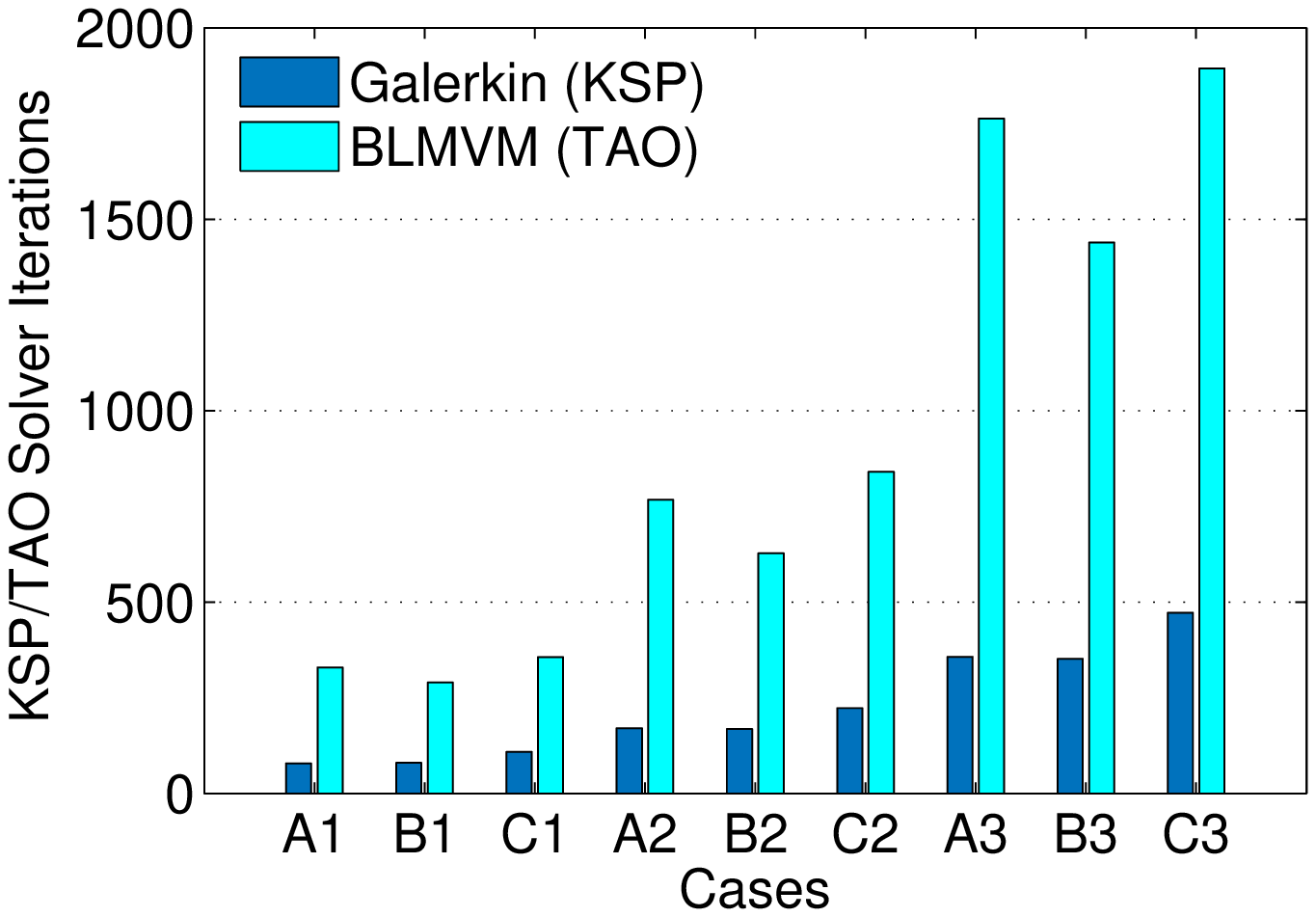}}
\caption{Cube with a hole: solver iterations needed for Galerkin and BLMVM.}
\label{Fig:S5_cube_iterations_galerkin}
\end{figure}
%------------------------------------------------;
%  Figure: Solver iterations - TRON              ;
%------------------------------------------------;
\begin{figure}[t]
\centering
\subfloat{\includegraphics[scale=0.55]{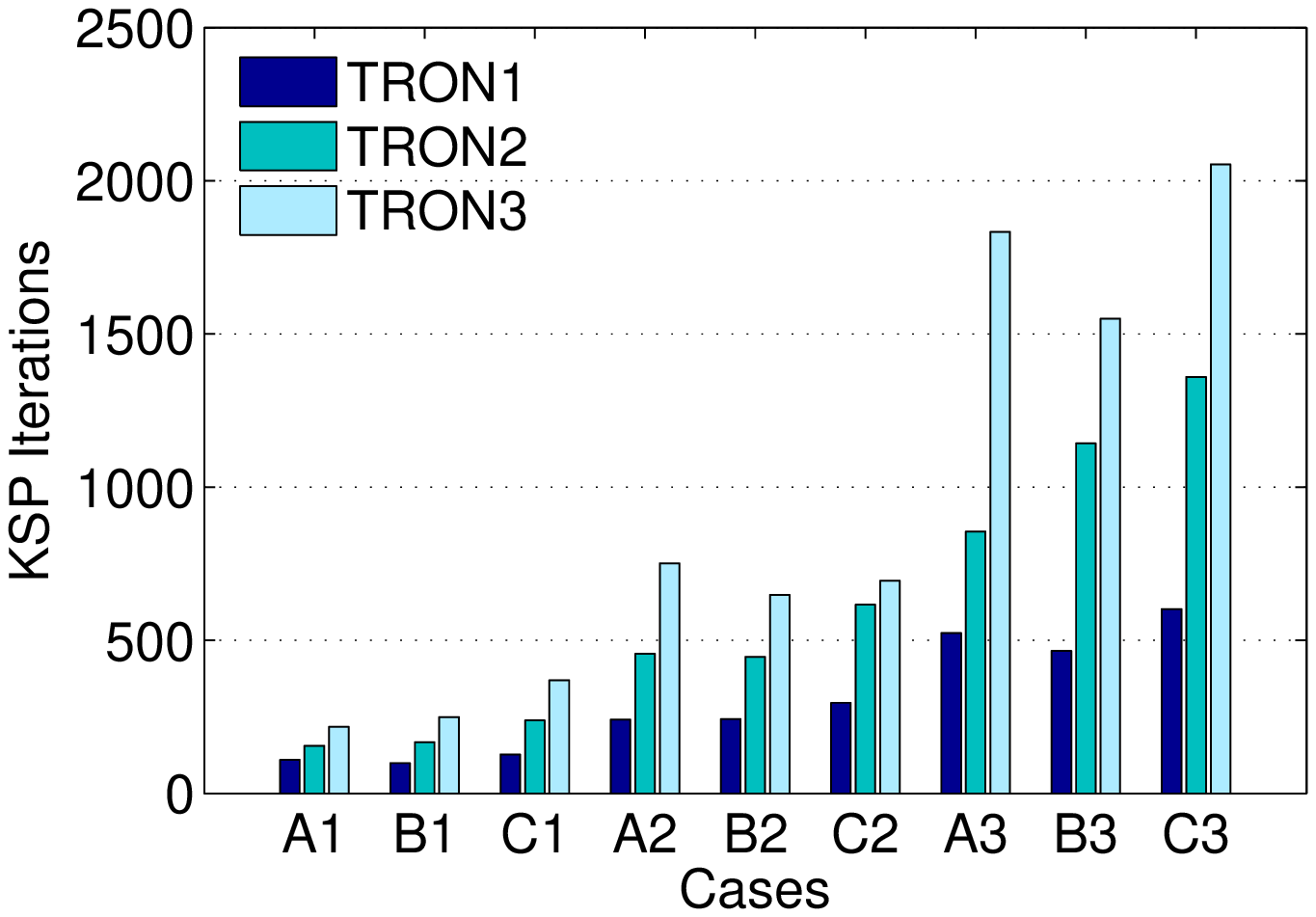}}
\subfloat{\includegraphics[scale=0.55]{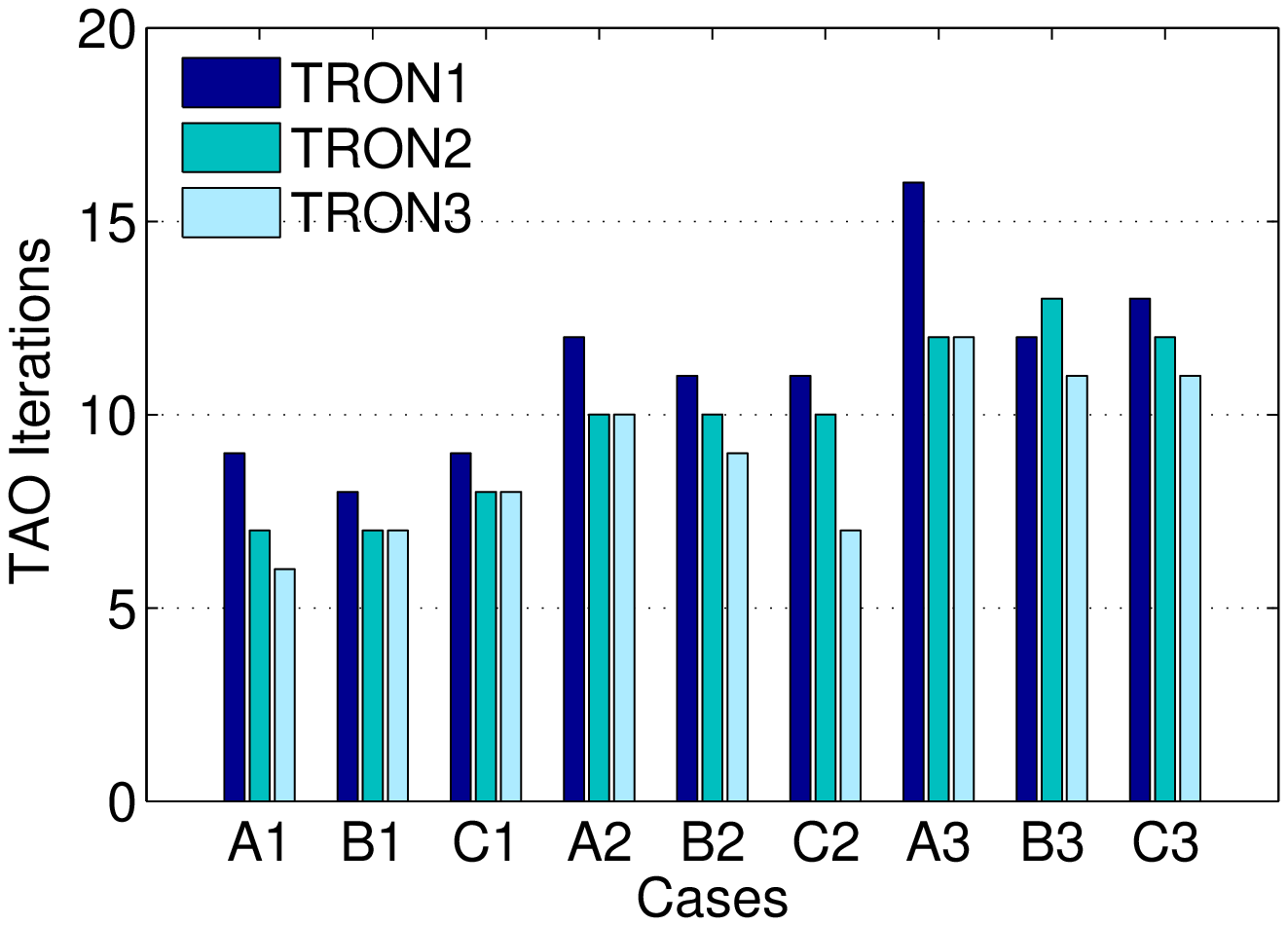}}
\caption{Cube with a hole: KSP (left) and TAO (right) solver iterations needed for TRON.}
\label{Fig:S5_cube_iterations_tron}
\end{figure}

%====================================================;
%  Subsubsection: computational efficiency           ;
%====================================================;
\subsubsection{Performance modeling}
We first consider the wall-clock time spent in the solvers on a single core. 
Table \ref{Tab:S5_cube_times_mustang} and \ref{Tab:S5_cube_times_wolf} depict 
the solver time for each mesh, and we first note that Mustang system requires 
significantly more wall-clock time to obtain a solution than Wolf; 
this behavior is expected due to the difference in HPC hardware specifications 
listed in Table \ref{Tab:S3_HPC}, specifically, Mustang has the lower clock rate 
and lower bandwidth (as determined through STREAMS Triad). It can also be seen that 
the various non-negative solvers consume varying amounts of wall-clock time.
BLMVM can require as much as ten times the amount of wall-clock time as the standard
Galerkin method. TRON on the other hand, does not consume nearly as much time but 
tightening the KSP tolerances will gradually increase the amount of time. 
We are interested in determine why these optimization solvers consume more 
wall-clock time, whether it be mostly due to additional workload 
associated with optimization-based techniques or due to the presence of 
relatively more complicated and expensive data structures compared to the 
standard solvers used for the Galerkin formulation. The first step is noting 
the total KSP and TAO iterations needed and how they vary with respect to problem size. 
Figure \ref{Fig:S5_cube_iterations_galerkin} depicts the KSP and TAO iterations for 
the Galerkin and BLMVM methods respectively. It is well-known that block 
Jacobi (also known as ILU(0)) requires more iterates as the size of the 
problem increases. In other words, the solver may exhibit poor scaling for 
extremely large problems, but we see that the BLMVM algorithm has an even poorer 
scaling rate of the solver iterates. For the TRON solvers, we document both 
the KSP and TAO iterates as shown in Figure \ref{Fig:S5_cube_iterations_tron}. 
We see that tightening the KSP tolerance increases the number of KSP iterates 
but requires slightly fewer TAO iterates. This behavior indicates that the more 
accurate the computed gradient projection is, the fewer 
optimization loops the solver has to perform.

%------------------------------------------------;
%  Figure: FLOPSS                                ;
%------------------------------------------------;
\begin{figure}[tp]
\centering
\subfloat[Mustang]{\includegraphics[scale=0.55]
{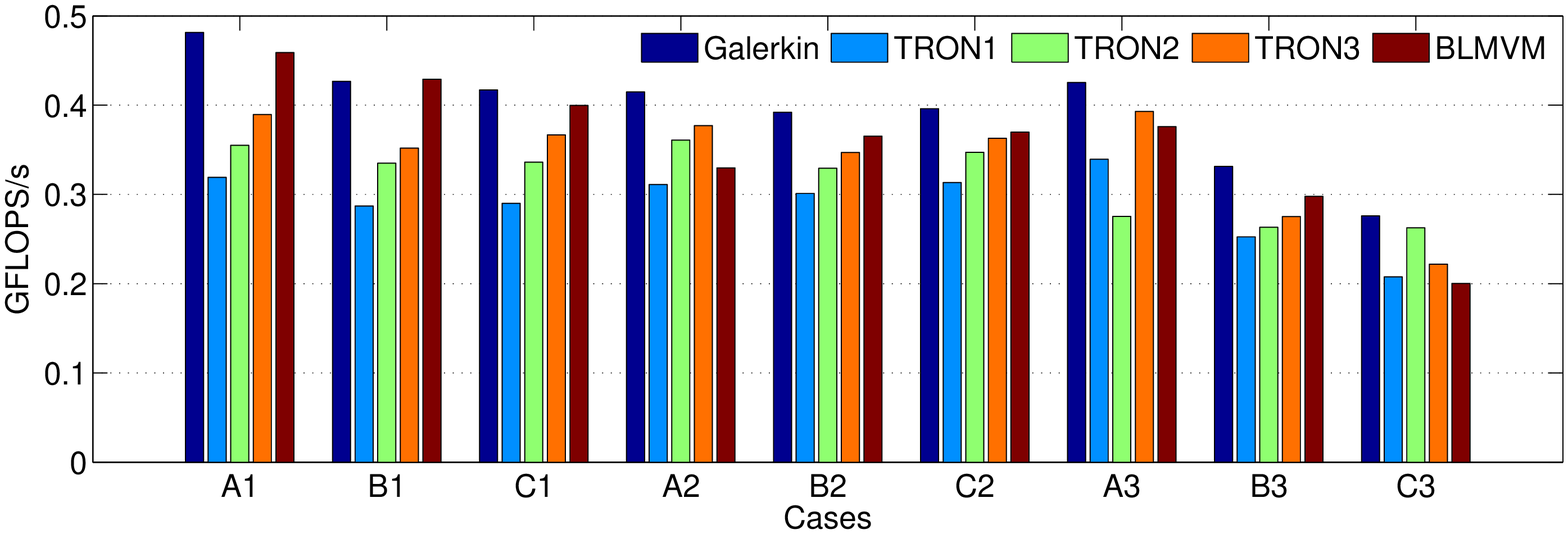}}
\label{Fig:S5_cube_flopss_mustang}
\subfloat[Wolf]{\includegraphics[scale=0.55]
{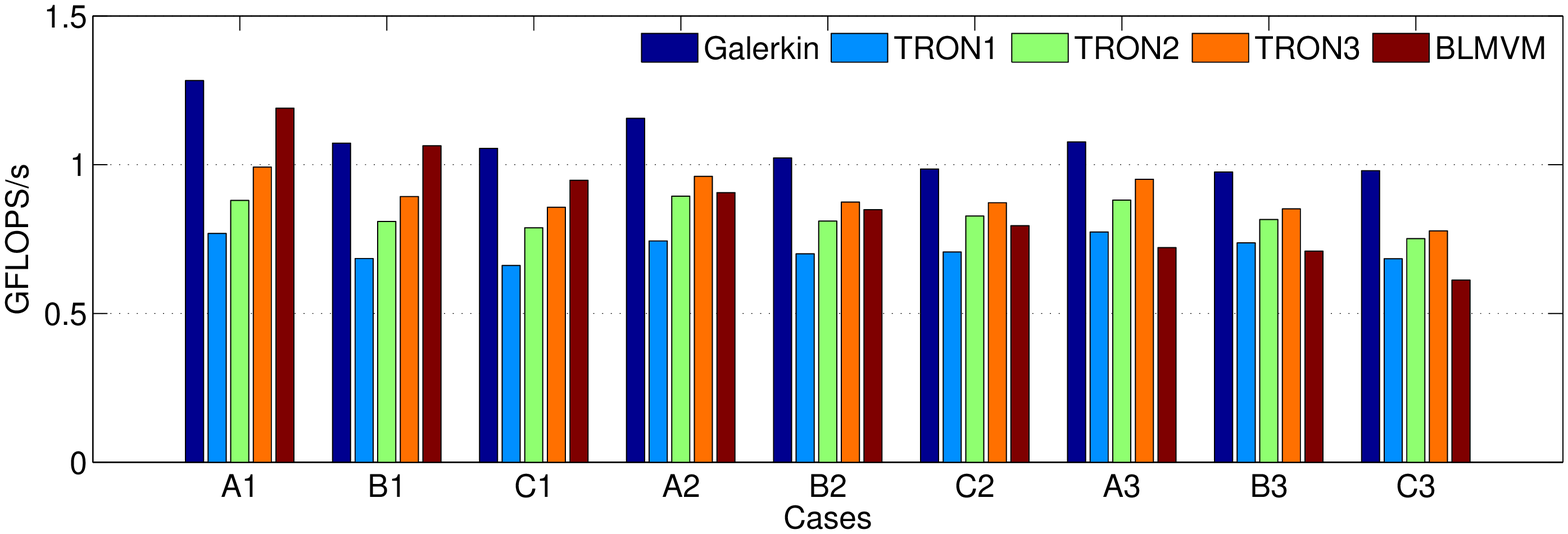}}
\label{Fig:S5_cube_flopss_wolf}\\
\caption{Cube with a hole: measured floating-point rate (FLOPS/s) on a single core.}
\label{Fig:S5_cube_flopss}
\end{figure}
%------------------------------------------------;
%  Figure: Arithmetic intensity for Cube         ;
%------------------------------------------------;
\begin{figure}[t]
\centering
\subfloat{\includegraphics[scale=0.55]{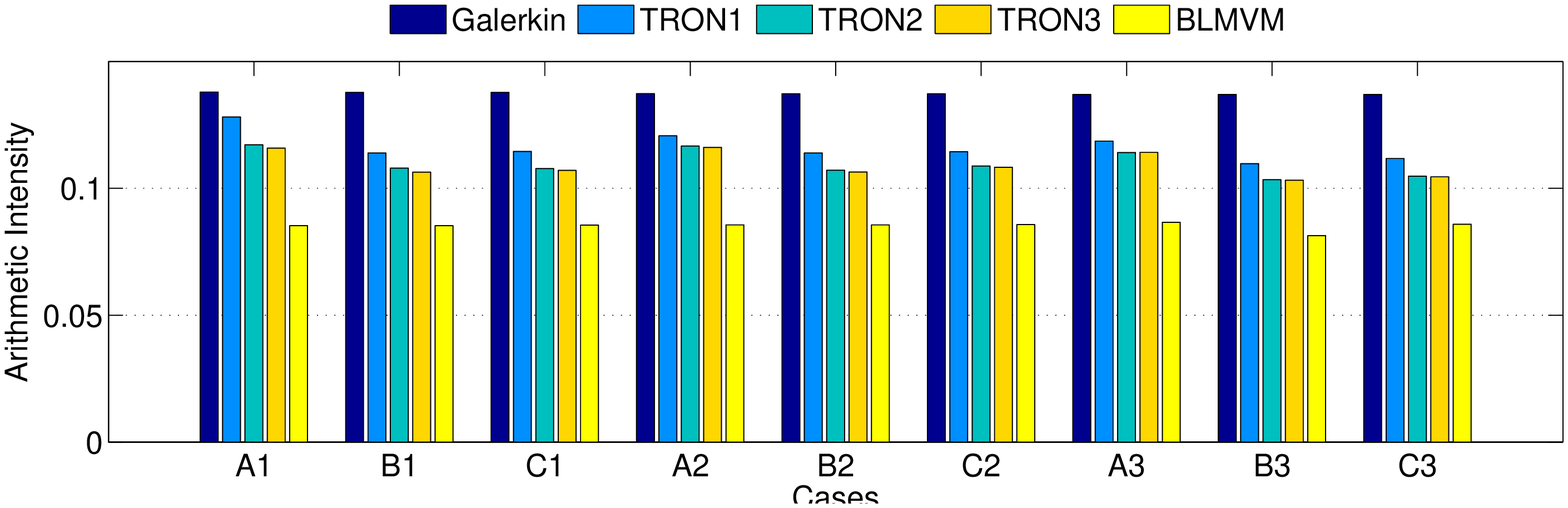}}
\caption{Cube with a hole: arithmetic intensity for all solvers and all cases  on a single processor.}
\label{Fig:S5_cube_ai}
\end{figure}

\begin{figure}[tp]
\centering
\subfloat[Mustang]{\includegraphics[scale=0.55]
{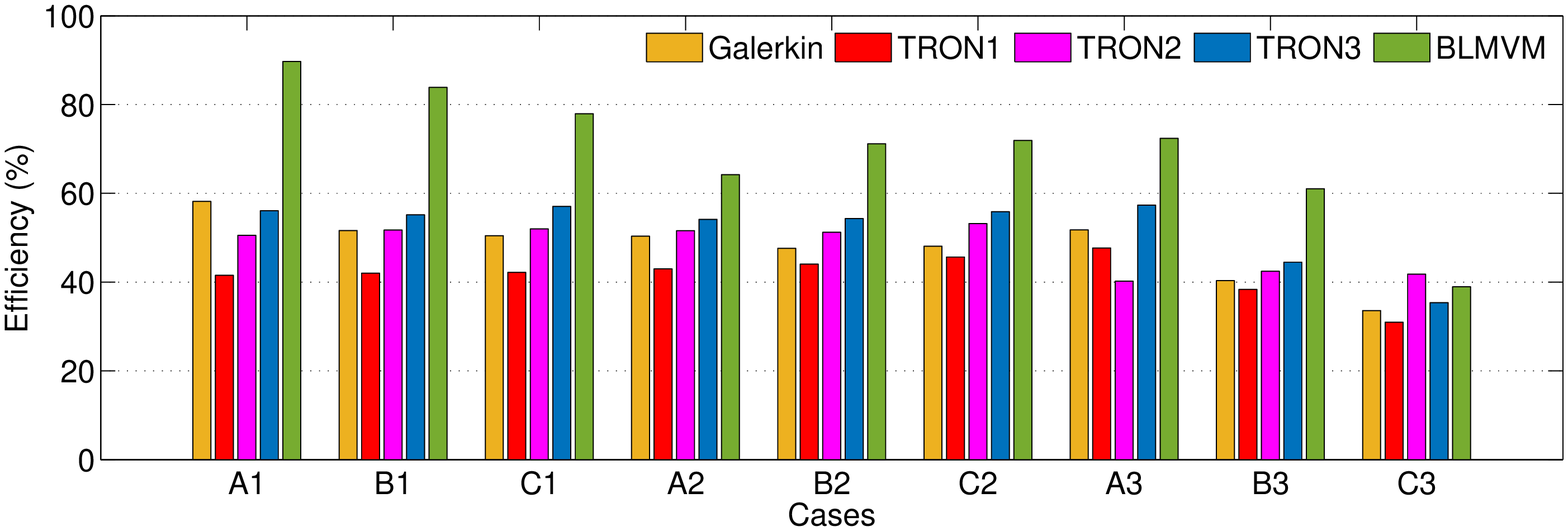}}
\label{Fig:S5_cube_flopss_mustang}
\subfloat[Wolf]{\includegraphics[scale=0.55]
{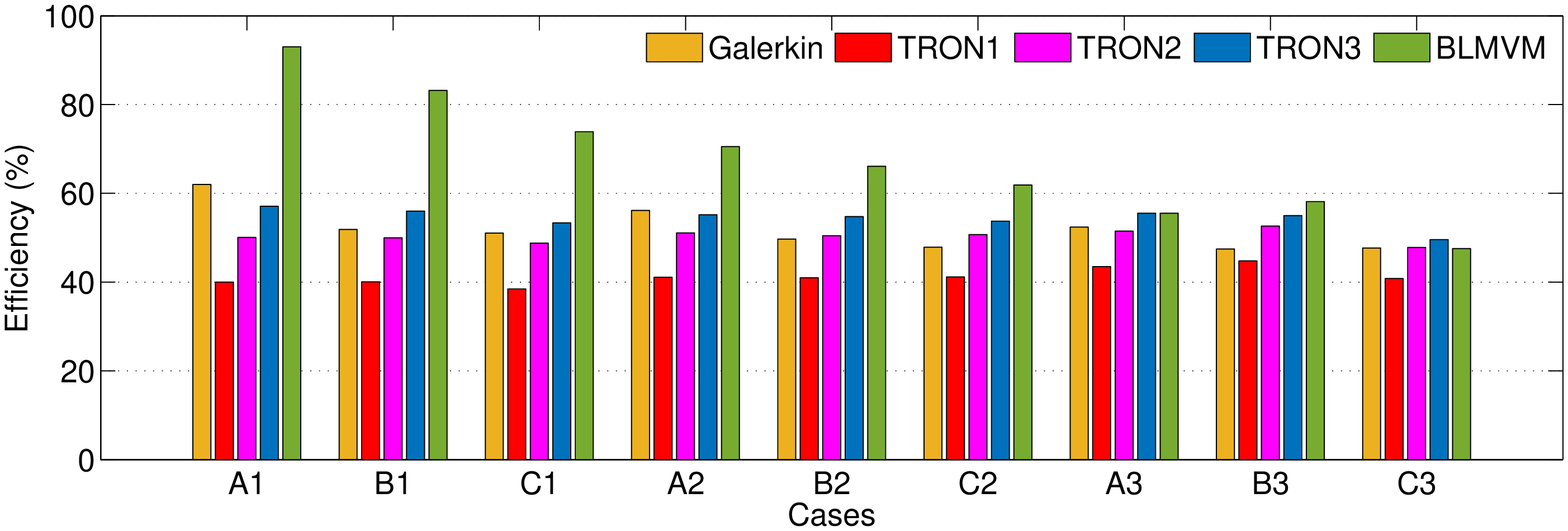}}
\label{Fig:S5_cube_flopss_wolf}\\
\caption{Cube with a hole: estimated floating-point efficiency with respect to the 
arithmetic intensity and measured memory bandwidth from STREAMS.}
\label{Fig:S5_cube_roofline}
\end{figure}

We also examine the measured floating-point rate provided by the PETSc performance logs,
as shown in Figure \ref{Fig:S5_cube_flopss}, of all five solvers across their respective 
machines, and the floating point performance decreases as the problem size grows. One 
could compare these numbers to the TPP and see that the hardware efficiencies are no 
greater than $5\%$, but it is difficult to draw any other conclusions with regard to the
computational performance. The calculated AI, based on our proposed performance model, 
is shown in Figure \ref{Fig:S5_cube_ai}. It is 
interesting to note that the AI remains largely invariant with problem size 
unlike the wall-clock time, solver iterations, and floating point rates. 
According to the perfect cache model, the Galerkin formulation's AI is greater 
than any of the non-negative methodologies. The optimization-based 
algorithm based on TAO's BLMVM solver has significantly more streaming/vector 
operations which explains the relatively lower AI. 
Using these metrics in equation \eqref{Eqn:S3_efficiency} as well as the 
STREAMS Triad bandwidth of one core as shown in Figure \ref{Fig:S3_streams}, 
the estimated roofline-based efficiencies are shown in Figure \ref{Fig:S5_cube_roofline}. 
Although the raw floating-point rate of BLMVM is lower than the Galerkin method, the roofline
model suggests that BLMVM is actually more efficient in the hardware sense. The TRON methods
have much lower floating-point rates, but these metrics can be improved or ``gamed" by 
tightening the KSP tolerances. This behavior leads us to believe that there is some 
latency associated with setting up the data structures needed to compute gradient 
descent projections.

%------------------------------------------------;
%  Figure: Efficiency                            ;
%------------------------------------------------;
%------------------------------------------------;
%  Figure: Speedup Mustang                       ;
%------------------------------------------------;
\begin{figure}[t]
\centering
\subfloat{\includegraphics[scale=0.55]{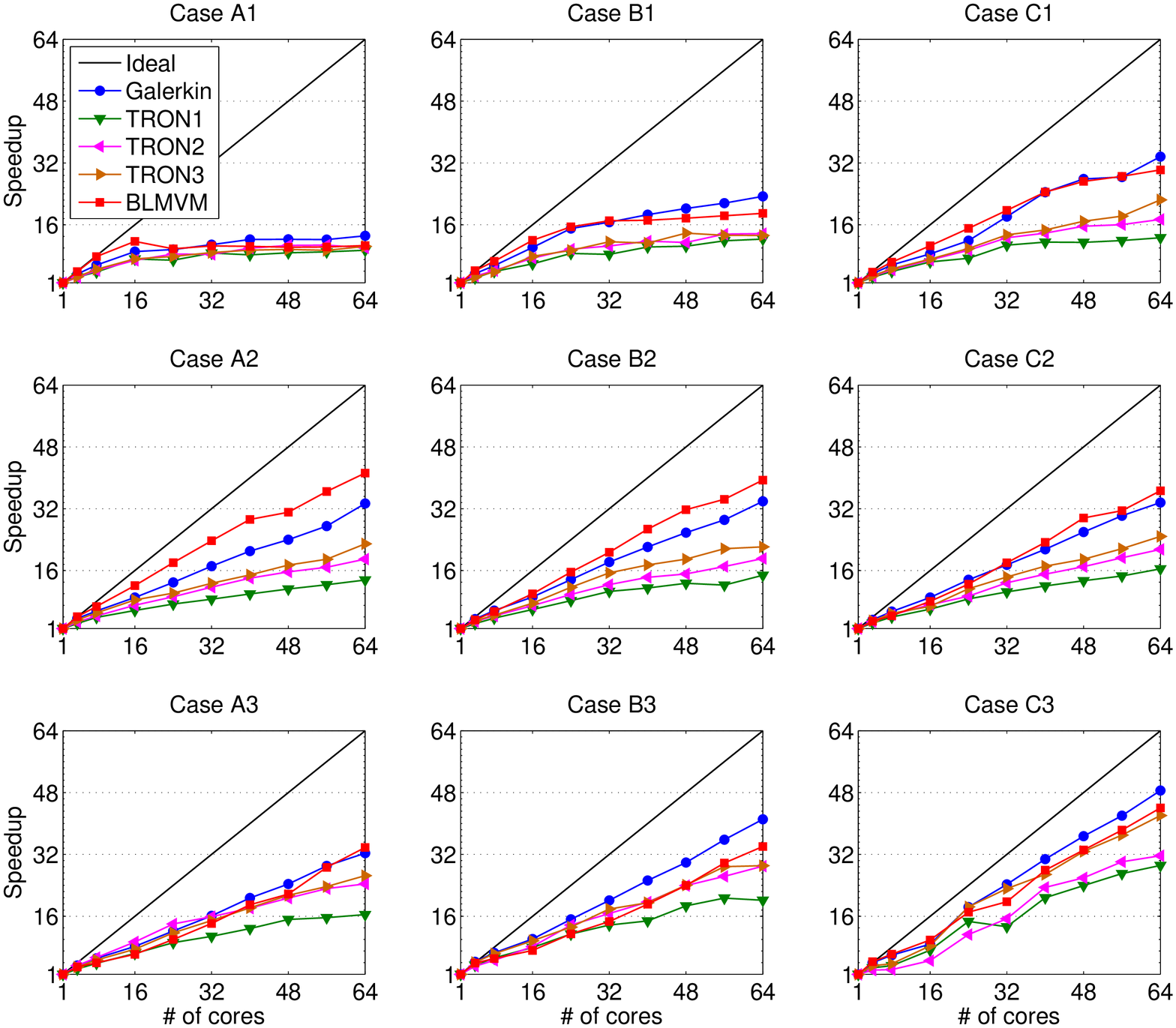}}
\caption{Cube with a hole: speedup for all 9 mesh cases up to 64 processors on the Mustang system (16 cores per node).}
\label{Fig:S5_cube_speedup_mustang}
\end{figure}
%------------------------------------------------;
%  Figure: Speedup Wolf                          ;
%------------------------------------------------;
\begin{figure}[t]
\centering
\subfloat{\includegraphics[scale=0.55]{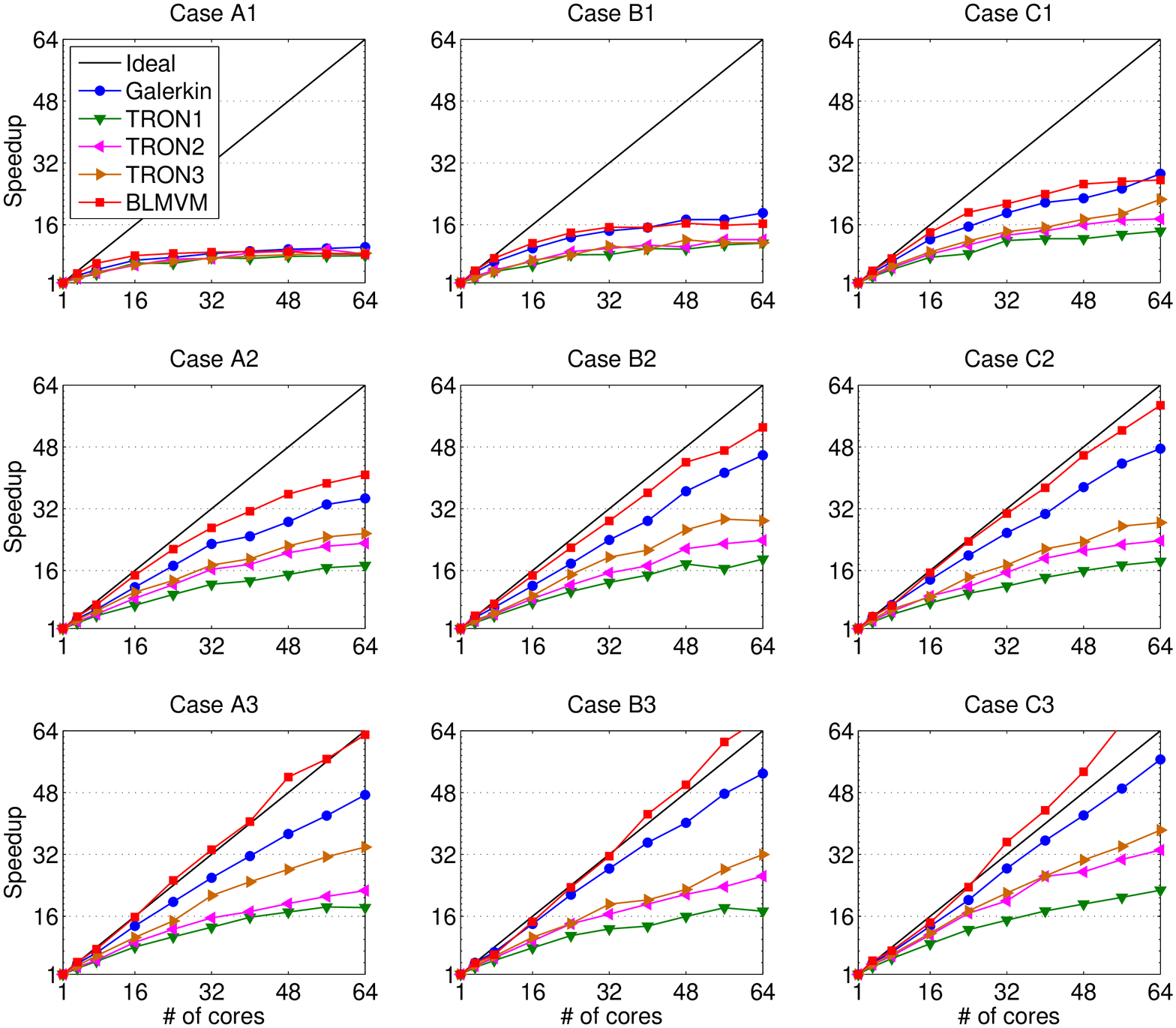}}
\caption{Cube with a hole: speedup for all 9 mesh cases up to 64 cores on the Wolf system (8 cores per node).}
\label{Fig:S5_cube_speedup_wolf}
\end{figure}
%====================================================;
%  Subsubsection: Strong-scaling                     ;
%====================================================;
\subsubsection{Strong-scaling}
The  metric of most interest to many computational scientists is the strong-scaling 
potential of any numerical framework. We conduct strong-scaling 
studies to measure the speedup of all nine case studies over 64 cores. 
Four Mustang nodes with 16 cores each and 8 Wolf nodes with 8 cores 
each are allocated for this study. We do not fully saturate the compute
nodes because the STREAMS benchmark indicates that there is little or no 
gain in memory performance when using a full node. Figure \ref{Fig:S5_cube_speedup_mustang} 
depicts the speedup on the Mustang system, and Figure \ref{Fig:S5_cube_speedup_wolf}
depicts the speedup on the Wolf system. First, we note that the parallel efficiency 
(actual speedup over ideal speedup) increases with problem size due to 
Amdahl's Law. We also note that Wolf exhibits better strong-scaling due to the
faster speedups for the same test studies. For all problems and machines, 
the TRON simulations are slightly less efficient in the parallel sense but can be 
improved by tightening the KSP tolerances. Interestingly, the BLMVM algorithm not 
only has the best roofline efficiency but also the best parallel speedup. 
We can infer from these results that although BLMVM is the more efficient 
optimization in the hardware sense, TRON is more efficient in the algorithmic sense 
due to its lesser time-to-solution. Our study has shown that one can draw 
correlations between the performance models conducted on a single-core and the 
actual speedup across multiple distributed memory nodes. As future solvers and 
algorithms are implemented within PETSc, we can use this performance model to 
assess how efficient they are in both the hardware and algorithmic sense and how 
efficiently they will scale in a parallel setting.

%------------------------------------------------;
%  Figure: Description of chromium subsurface    ;
%------------------------------------------------;
\begin{figure}[t]
\centering
\subfloat{\includegraphics[scale=0.45]{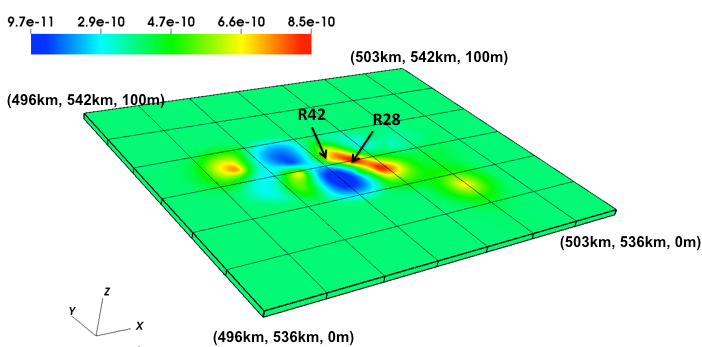}}
\caption{Chromium plume migration in the subsurface: Permeability field (m$^2$) and the locations of the pumping well (R28) and contaminant source (R42).}
\label{Fig:S5_chromium_description}
\end{figure}
%===================================================;
%  Subsection: Transport of chromium in subsurface  ;
%===================================================;
\subsection{Transport of chromium in subsurface}
%------------------------------------------------;
%  Table: Chromium properties                    ;
%------------------------------------------------;
\begin{table}[b]
  \centering
  \caption{Chromium plume migration in the subsurface: parameters \label{Tab:S5_chrom_properties}}
  \begin{tabular}{ll}
    \hline
    Parameter & Value \\
    \hline
    $\alpha_L$ & 100 m \\
    $\alpha_T$ & 0.1 m \\
    Contaminant source (R42) & $1\times 10^{-4}$ kg/m$^2$s$^2$ \\
    $\Delta t$ & 0.2 days \\
    Domain size & 7000 km$\times$6000 km$\times$100 m \\
    $D_M$ & $1\times 10^{-9}$ m$^{2}$/s \\
    Permeability & Varies \\
    Pumping well (R28) & -0.01 kg/m$^2$s$^2$\\
    Total hexahedrons & 1,984,512 \\
    Total vertices & 2,487,765 \\
    $\mathbf{v}$ & Varies with position  \\
    Viscosity & 3.95$\times 10^{-5}$ Pa s \\
    \hline
  \end{tabular}
\end{table}
%------------------------------------------------;
%  Figure: Numerical results chromium orig       ;
%------------------------------------------------;
\begin{figure}[htp]
\centering
\subfloat{\includegraphics[scale=0.36]{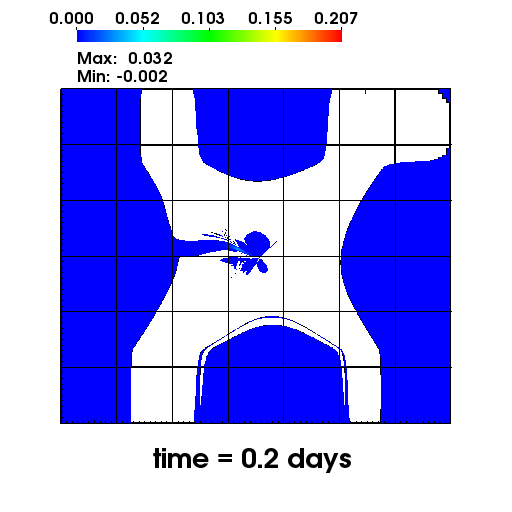}}
\subfloat{\includegraphics[scale=0.36]{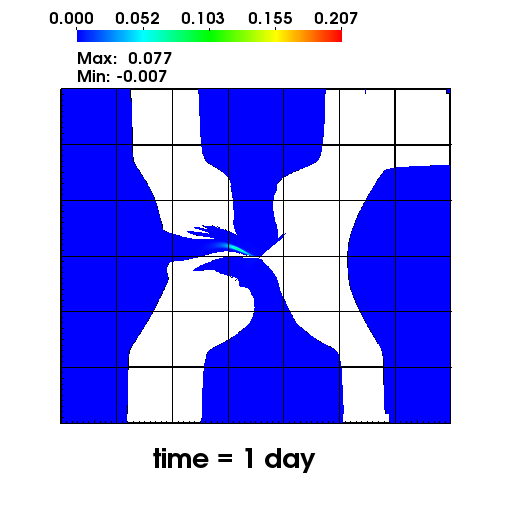}}\\
\subfloat{\includegraphics[scale=0.36]{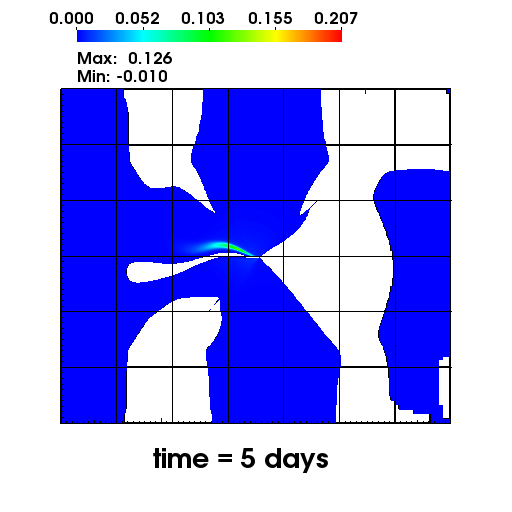}}
\subfloat{\includegraphics[scale=0.36]{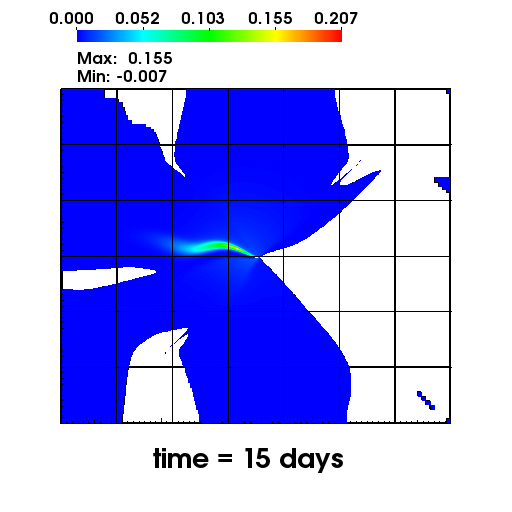}}\\
\subfloat{\includegraphics[scale=0.36]{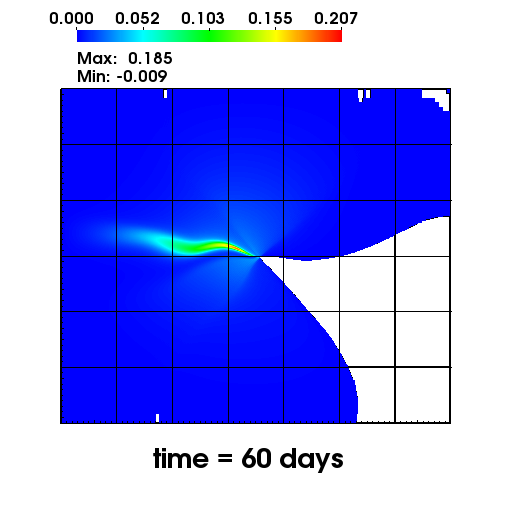}}
\subfloat{\includegraphics[scale=0.36]{Figures/chrom_orig_900.png}}
\caption{Chromium plume migration in the subsurface: concentrations at select times using the Galerkin formulation.}
\label{Fig:S5_chromium_orig}
\end{figure}
%------------------------------------------------;
%  Figure: Numerical results chromium nonneg     ;
%------------------------------------------------;
\begin{figure}[htp]
\centering
\subfloat{\includegraphics[scale=0.36]{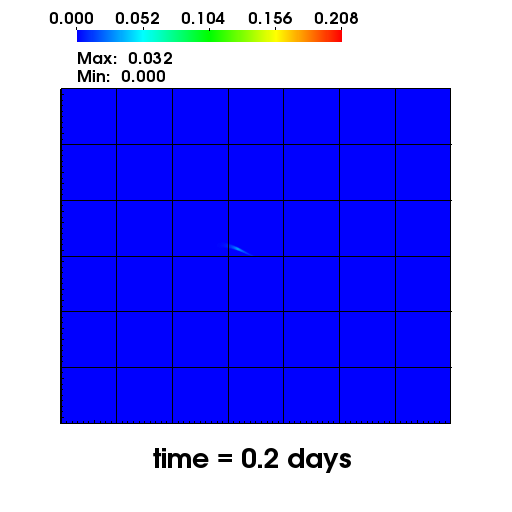}}
\subfloat{\includegraphics[scale=0.36]{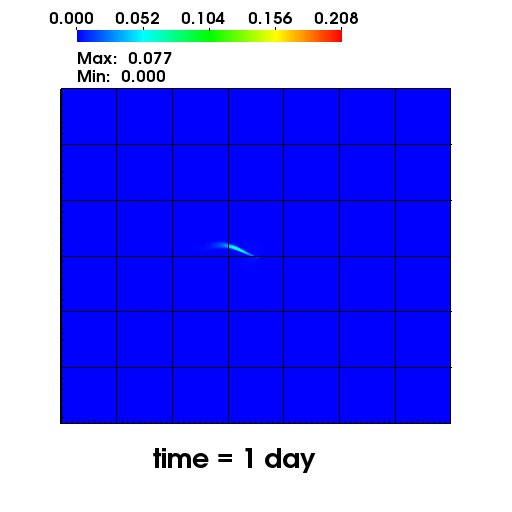}}\\
\subfloat{\includegraphics[scale=0.36]{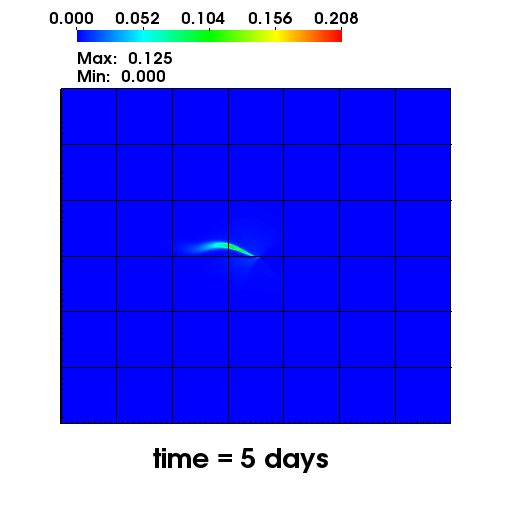}}
\subfloat{\includegraphics[scale=0.36]{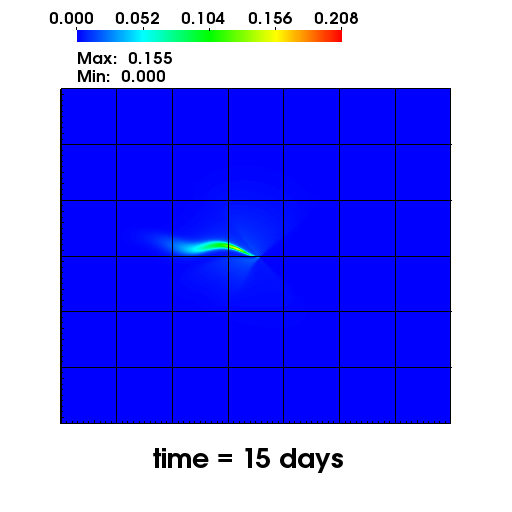}}\\
\subfloat{\includegraphics[scale=0.36]{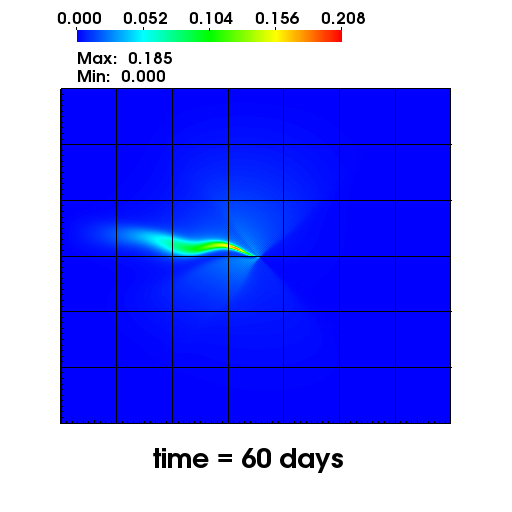}}
\subfloat{\includegraphics[scale=0.36]{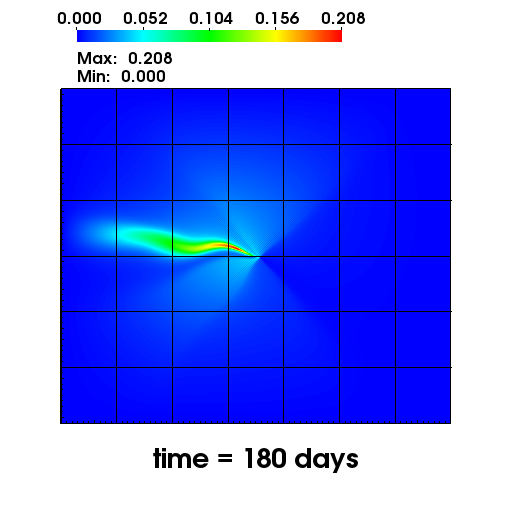}}
\caption{Chromium plume migration in the subsurface: concentrations at select times using the non-negative methodology.}
\label{Fig:S5_chromium_nonneg}
\end{figure}
Subsurface clean-up due to anthropogenic contamination 
is a big challenge \cite{us2004cleaning}.  
Remediation studies \cite{hammond2010field,harp2013contaminant} need accurate predictions 
of transport of the involved chemical species, which are obtained 
using limited data at monitoring wells and through 
numerical simulations.
% Therefore, an accurate 
%predictive numerical simulator that predicts the 
%fate of a chemical species in the subsurface is 
%of significant importance. 
%
To accurately predict the fate of the contaminant,
a transport solver that: a) is robust, in that it will not 
give unphysical solutions, and b) can handle field-scale 
scenarios, is needed. The computational framework that is 
proposed in this paper is an ideal candidate 
for such problems. We now consider a realistic 
large-scale problem to predict the fate of 
chromium in the Los Alamos, New Mexico area. 
The chromium was released into the Sandia canyon in the 50s up to 
early 70s. Back then chromium was used as an anti-corrosion
agent for the cooling towers at a power plant at the 
Los Alamos National Laboratory (see \cite{heikoop2014isotopic} 
and references therein for details).   

Here we study the effectiveness of our proposed framework
to this real-world scenario of predicting the extent of 
chromium plume. The following is a conceptual model domain
that is considered:  A domain of size
 $[496,503]\mathrm{km} \times [536,542]
\mathrm{km} \times [0,100]\mathrm{m}$ with the permeability field (m$^2$) as shown 
in Figure \ref{Fig:S5_chromium_description}. R42 in 
Figure \ref{Fig:S5_chromium_description} is estimated to be the contaminant source
location and a pumping well is located at R28. The parameters used for this 
problem are shown in Table \ref{Tab:S5_chrom_properties}, and we employ the 
following boundary conditions:
%-------------------------------------------------;
%  Equation: boundary conditions for Chromium     ;
%-------------------------------------------------;
\begin{align}
 \label{Eqn:S5_chrom_bc_conc}
 &c^{\mathrm{p}}(x=496\mathrm{km},y,z) = c^{\mathrm{p}}(x=503\mathrm{km},y,z) = 
c^{\mathrm{p}}(x,y=536\mathrm{km},z) =  c^{\mathrm{p}}(x,y=542\mathrm{km},z) = 0
\end{align}
For this highly heterogeneous problem, we 
employ PETSc's algebraic multi-grid preconditioner (GAMG) and
couple this with the TRON algorithm for the non-negative solver. Our goal is 
to examine its strong-scaling potential across 1024 cores. 
We first solve the steady flow equation (based on mass balance and
Darcy's model to relate pressure and mass flux) with the pumping well 
located at R28. Cell-wise velocity is obtained from the resulting pressure field and 
used to calculated element-wise dispersion tensor.  We then solve 
the transient diffusion problem (with tensorial dispersion) with a 
constant contaminant source located at R40 for up to 180 days. 
The concentrations at select time levels for Galerkin formulation and 
non-negative methodology are shown in Figures \ref{Fig:S5_chromium_orig} 
and \ref{Fig:S5_chromium_nonneg}, respectively. Negative concentrations arise with 
the Galerkin formulation even as the solution approaches steady-state.

\begin{figure}[t]
\centering
\subfloat{\includegraphics[scale=0.55]{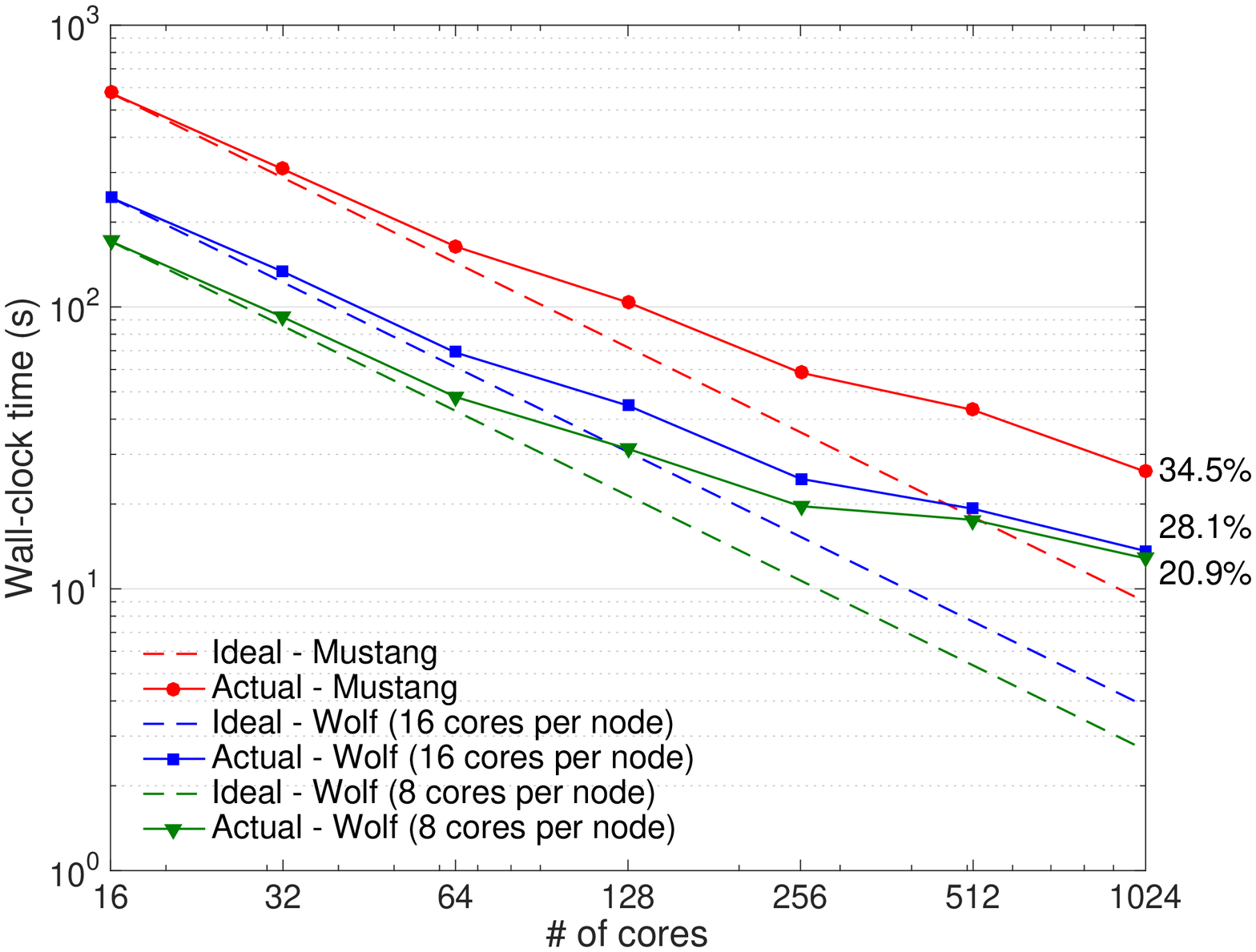}}
\caption{Chromium plume migration in the subsurface: wall-clock time of the TRON 
optimization solver with multi-grid preconditioner (GAMG) versus number 
of processors after the first time level. Two Wolf cases are considered, where 
we fully saturate a compute node (16 cores) and where we partially saturate a 
compute node (8 cores per node). The parallel efficiency with respect to 16 
cores is shown on the right hand side.}
\label{Fig:S5_chromium_strongscale_first}
\end{figure}

Figure \ref{Fig:S5_chromium_strongscale_first} depicts the amount of wall-clock 
time with respect to the number of cores at the first time level. We see here 
that the demonstrates good parallel performance across 1024 cores with up to 
35 percent parallel efficiency. Unlike the previous benchmark problem, we consider 
a case where we completely saturate a Wolf node by using all 16 cores and notice
that the performance is slightly worse than using a partially saturated node 
(8 cores). This change of behavior can be attributed to the lack of memory 
performance improvement one achieves when using all 16 cores as shown in 
Figure \ref{Fig:S3_streams}. Interprocess communication becomes a major component of
the Wolf simulations so the parallel scalability gets worse the more efficiently TRON
conducts its work. Nonetheless, the higher quality computing resources 
of a Wolf node results in faster solve times than the solve time on Mustang 
even with lesser parallel efficiency. 

\begin{figure}[t]
\centering
\subfloat{\includegraphics[scale=0.55]{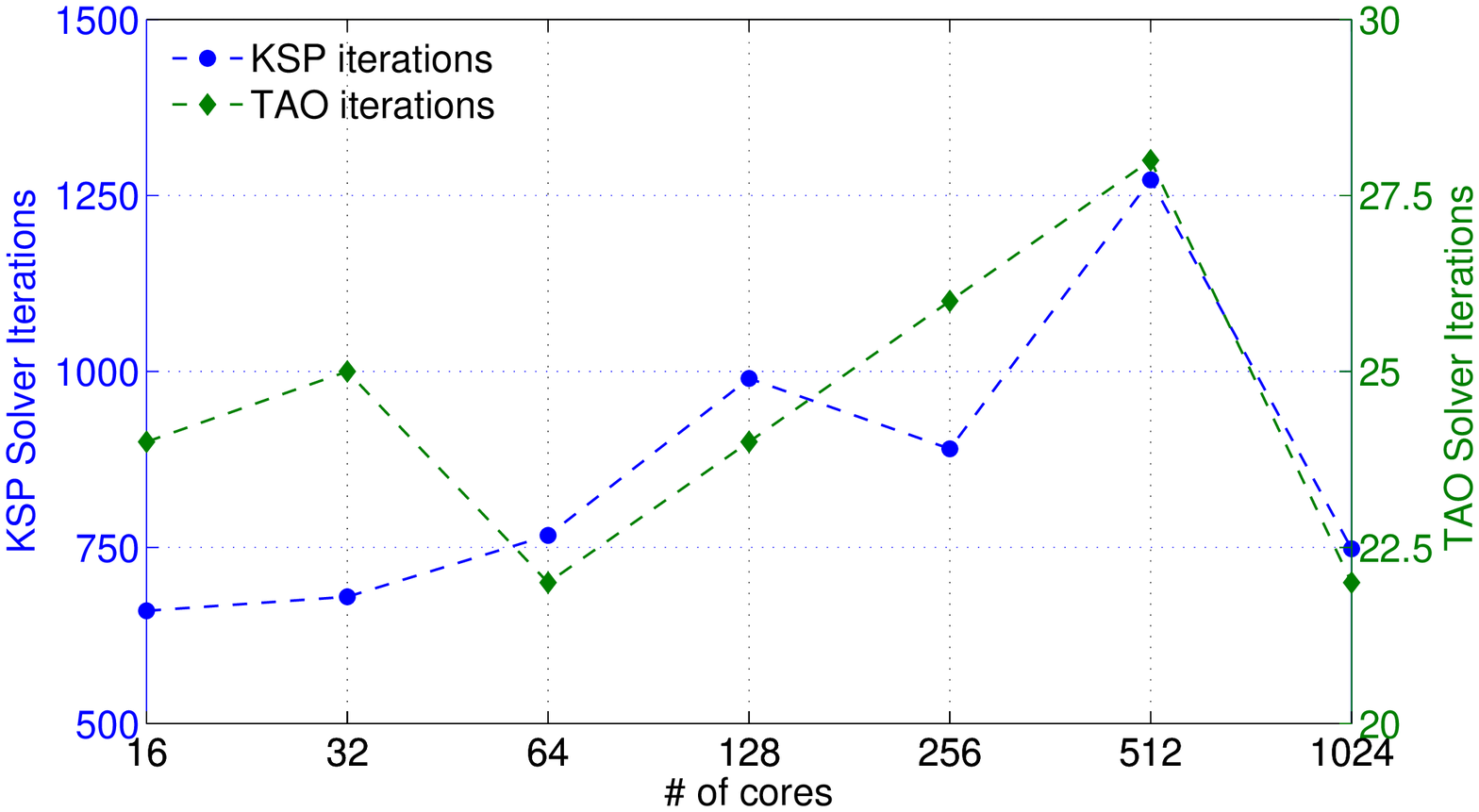}}
\caption{Chromium plume migration in the subsurface: number of KSP (left hand side) 
and TAO (right hand side) solver iterations for the TRON optimization solver versus 
number of cores after the first time level.}
\label{Fig:S5_chromium_iterations_first}
\end{figure}

Another metric of interest is the number of solver iterations required for
convergence across various number of MPI processes. Figure 
\ref{Fig:S5_chromium_iterations_first} depicts the number of KSP solver iterations
and TAO solver iterations across 1024 cores, and we notice that there are 
significant fluctuations. This trend is largely attributed to the accumulation 
of numerical round-offs from the TRON algorithm. One can reduce
these fluctuations by tightening the solver tolerances, but the strong-scaling 
remains largely unaffected even for the results shown. This study suggests that 
the proposed non-negative methodology using TRON with GAMG preconditioning 
is suitable for large-scale transient heterogeneious and anisotropic diffusion 
equations.

%------------------------------------------------;
%  Figure: Chromium strong scale wolf first      ;
%------------------------------------------------;

%------------------------------------------------;
%  Figure: Chromium strong scale wolf first      ;
%------------------------------------------------;

%*********************************************;
%                                             ;
%  NAME                                       ;
%    S5_LargeScale_RM.tex          	      ;
%                                             ;
%  WRITTEN BY                                 ;
%	 Justin Chang			      ;
%	 Satish Karra		              ;
%    Kalyana Babu Nakshatrala                 ;
%                                             ;
%*********************************************;
\section{CONCLUDING REMARKS}
\label{Sec:Concluding_Remarks}
We presented a parallel non-negative computational framework 
suitable for solving large-scale steady-state and transient 
anisotropic diffusion equations. The main contribution is 
that the proposed parallel computational framework satisfies 
the discrete maximum principles for large-scale diffusion-type 
equations even on general computational grids. 
The parallel framework is built upon PETSc's DMPlex
data structure, which can handle unstructured meshes, 
and TAO for solving the resulting optimization problems 
from the discretization formulation. We have conducted 
systematic performance modeling and strong-scaling 
studies to demonstrate the efficiency, both in the parallel and
hardware sense of the computational framework. The robustness of 
the proposed framework has been illustrated by solving a 
large-scale realistic problem involving the transport of 
chromium in the subsurface at Los Alamos, New Mexico.
Future areas of research include: 
(a) extending the proposed parallel framework to 
advective-diffusive and advective-diffusive-reactive 
systems, and (b) posing the discrete problem as a 
variational inequality, which will be valid even 
for non-self-adjoint operators, and use other PETSc
capabilities to solve such variational 
inequalities.

%%============================;
%%  Section: Acknowledgments  ;
%%============================;
\section*{ACKNOWLEDGMENTS}
The authors thank Matthew G. Knepley (Rice University) 
for his invaluable advice. The authors also thank
the Los Alamos National Laboratory (LANL)
Institutional Computing program. 
JC and KBN acknowledge the financial support from the 
Houston Endowment Fund and from the Department of Energy 
through Subsurface Biogeochemical Research Program. 
SK thanks the LANL LDRD program and the LANL Environmental
Programs Directorate for their support.
The opinions expressed in this paper are those of 
the authors and do not necessarily reflect that 
of the sponsors.

%================;
%  Bibliography  ;
%================;
\bibliographystyle{plain}
\bibliography{bibliography,Master_References/Books}
\newpage
\clearpage
%================================;
%  Include all the figures here  ;
%================================;
%%
\end{document}